\documentclass[a4paper,10pt]{article}
\usepackage{stmaryrd}
\usepackage{amsfonts}
\usepackage{bbm}
\usepackage{amscd}
\usepackage{mathrsfs}
\usepackage{latexsym,amssymb,amsmath,amscd,amscd,amsthm,amsxtra}
\usepackage[dvips]{graphicx}
\usepackage[utf8]{inputenc}
\usepackage[T1]{fontenc}
\usepackage{lmodern}
\usepackage{amssymb}
\usepackage[all]{xy}
\usepackage{nicefrac,mathtools,enumitem}
\usepackage{microtype}

\newtheorem{thm}{Theorem}[section]
\newtheorem{lem}[thm]{Lemma}
\newtheorem{cor}[thm]{Corollary}
\newtheorem{pro}[thm]{Proposition}

\newtheorem{rmk}[thm]{Remark}

\newtheorem{defi}[thm]{Definition}

\newcommand{\lon }{\,\rightarrow\,}
\newcommand{\be }{\begin{equation}}
\newcommand{\ee }{\end{equation}}
\newcommand {\emptycomment}[1]{} 

\newcommand{\pf}{\noindent{\bf Proof.}\ }

\newcommand{\huaB}{\mathcal{B}}

\newcommand{\frkad}{\mathfrak {ad}}

\newcommand{\huaH}{\mathcal{H}}

\newcommand{\huaZ}{\mathcal{Z}}

\newcommand{\frkg}{\mathfrak g}\newcommand{\g}{\mathfrak g}
\newcommand{\frkh}{\mathfrak h}

\newcommand{\frkk}{\mathfrak k}

\newcommand{\frkCM}{\mathfrak {CM}}
\newcommand{\frkSL}{\mathfrak {SL}}

\def\qed{\hfill ~\vrule height6pt width6pt depth0pt}

\newcommand{\Id}{\rm{Id}}

\newcommand{\im}{\mathrm{\im}}

\newcommand{\dM}{\mathrm{d}}
\newcommand{\E}{\mathrm{E}}

\newcommand{\Hom}{\mathrm{Hom}}
\newcommand{\SDer}{\mathrm{SDer}}
\newcommand{\SInn}{\mathrm{SInn}}
\newcommand{\cen}{\mathsf{Cen}}
\newcommand{\Der}{\mathrm{Der}}
\newcommand{\Ext}{\mathrm{Ext}}
\newcommand{\cm}{\mathrm{CM}}
\newcommand{\Sl}{\mathrm{SL}}
\newcommand{\A}{\mathrm{A}}
\newcommand{\OM}{\mathrm{O}}

\newcommand{\DER}{\mathrm{DER}}
\newcommand{\D}{\mathrm{D}}

\newcommand{\SOut}{\mathrm{SOut}}

\newcommand{\gl}{\mathfrak {gl}}

\newcommand{\End}{\mathrm{End}}
\newcommand{\ad}{\mathrm{ad}}
\newcommand{\pr}{\mathrm{pr}}

\newcommand{\Img}{\mathrm{Im}}

\newcommand{\h}{\frkh}

\begin{document}
\title{
{Cohomological characterizations of non-abelian extensions of strict Lie 2-algebras\thanks
 {
Research supported by NSFC (11471139) and NSF of Jilin Province (20170101050JC).
}}}

\author{ Rong Tang and Yunhe Sheng \\
Department of Mathematics, Jilin University,\\
 Changchun 130012, Jilin, China
\\\vspace{3mm}
Email: tangrong16@mails.jlu.edu.cn,~ shengyh@jlu.edu.cn}

\date{}


\maketitle

\begin{abstract}
In this paper, we study non-abelian extensions of strict Lie 2-algebras via the cohomology theory. A non-abelian extension of a strict Lie 2-algebra $\g$ by $\frkh$ gives rise to a strict homomorphism from $\g$ to $\SOut(\frkh)$. Conversely, we prove that the obstruction of existence of non-abelian
extensions of strict Lie 2-algebras associated to a strict Lie 2-algebra homomorphism from $\g$ to $\SOut(\frkh)$ is given by an element in the third cohomology group. We further
prove that the isomorphism classes of non-abelian extensions of   strict Lie 2-algebras  are classified by the second cohomology group.

\end{abstract}

\section{Introduction}
Eilenberg and Maclane \cite{EM} developed a theory of non-abelian extensions of abstract groups in the 1940s, leading to the low dimensional non-abelian group cohomology. Then there are a lot of analogous results for Lie algebras \cite{AMR,F,Hochschild,IKL}. Non-abelian extensions of Lie algebras can be described by some linear maps regarded as derivations of Lie algebras. This result was generalized to the case of super Lie algebras in \cite{AMR2}, to the case of topological Lie algebras in \cite{Hermann} and to the case of Lie algebroids in \cite{brahic}.

Lie 2-algebras are the categorification of Lie algebras \cite{BC}. In a Lie 2-algebra, the Jacobi identity is replaced by a natural isomorphism, which satisfies its own coherence law, called the Jacobiator identity. The 2-category of Lie 2-algebras is equivalent to the 2-category of 2-term $L_{\infty}$-algebras \cite{BC,LS,LM,stasheff:shla}, so people also view a 2-term $L_{\infty}$-algebras as a Lie 2-algebra directly.

Some work on the study of non-abelian extensions of Lie 2-algebras has been done recently. In \cite{CSZ}, the authors showed that a non-abelian extension of a Lie 2-algebra $\g$ by a Lie 2-algebra $\frkh$ can be characterized by a Lie 3-algebra homomorphism from $\g$ to the derivation Lie 3-algebra $\DER(\frkh)$. In \cite{lazarev:extension}, the author classified general non-abelian extensions of $L_\infty$-algebras via homotopy classes of $L_\infty$-morphism. On the other hand, strict Lie 2-algebras are in one-to-one correspondence with Lie algebra crossed modules \cite{BC}. Actor of Lie algebra crossed modules was constructed in \cite{CL}. Furthermore, the authors proved that associated a non-abelian extension of Lie algebra crossed modules $(\g_0,\g_1,\dM_{\g})$ by $(\frkh_0,\frkh_1,\dM_{\frkh})$, there is a Lie algebra crossed module morphism from $(\g_0,\g_1,\dM_{\g})$ to the out actor of $(\frkh_0,\frkh_1,\dM_{\frkh})$. However, whether there is obstruction of existence of a non-abelian extension associated to such a morphism and how to classify such extensions are still open problems.

The purpose of  this paper is to solve the above problems using the language of strict Lie 2-algebras.   To do that, first we outline the theory of strict derivations of strict Lie 2-algebras. We construct a strict Lie 2-algebra $\SDer(\frkg)$ associated to strict derivations, which plays important role when we consider non-abelian extensions of strict Lie 2-algebras. This part is not totally new, and one can obtain these results from \cite{CSZ,DL}. Then we show that a non-abelian extension of a strict Lie 2-algebra $\g$ by $\h$ naturally gives a strict homomorphism from $\g$ to $\SOut(\h)$. In the case of $\cen(\h)=0$, there is a one-to-one correspondence between isomorphism classes of non-abelian extensions of $\g$ by $\h$ and strict homomorphisms from  $\g$ to $\SOut(\h)$. Finally, we consider the obstruction of existence of non-abelian extensions of $\g$ by $\h$ associated to a strict homomorphism from  $\g$ to $\SOut(\h)$ in the case that $\cen(\g)\neq 0$. We show that the obstruction is given by an element in the third cohomology group. Furthermore, we use the second cohomology group to classify non-abelian extensions once they exist.

The paper is organized as follows. In Section 2, we give precise formulas for representations and cohomologies of strict Lie 2-algebras. In Section 3, we study strict derivations of strict Lie 2-algebras. In Section 4,
 by choosing a section, we give a general description of a non-abelian extension of a strict Lie 2-algebra $\g$ by $\frkh$. In Section 5, we prove that when $\cen(\frkh)=0$, there is a one-to-one correspondence between   isomorphism classes of non-abelian extensions of a strict Lie 2-algebra $\g$ by $\frkh$  and   strict Lie 2-algebra homomorphisms from $\g$ to the Lie 2-algebra $\SOut(\frkh)$. In Section 6, we identify a
third cohomological obstruction to the existence of non-abelian extensions associated to a homomorphism from $\g$ to the Lie 2-algebra $\SOut(\frkh)$. Furthermore,  we classify non-abelian extensions once they exist by the second cohomology group $\huaH^2(\g;\cen(\frkh))_{\hat{\mu}}$.

\section{Preliminary}
In this section we recall some basic concepts of representations and cohomologies of strict Lie 2-algebras.
\begin{defi}\label{slie2}
A strict Lie $2$-algebra  is a  $2$-term  graded vector spaces $ \g= \g_{1}\oplus
  \g_{0}$ equipped  with a linear map $\dM_\g:\g_1\longrightarrow \g_0$ and a skew-symmetric bilinear map $[\cdot,\cdot]_\g:\g_i\wedge\g_j\longrightarrow \g_{i+j},0\le i+j\le 1$, such that for all $x,y,z\in \g_0,~a,b\in \g_1$   the following equalities are satisfied:
\begin{itemize}
\item[\rm(a)] $\dM_{\g} [x,a]_\g=[x,\dM_{\g} a]_\g$,
\item[\rm(b)] $[\dM_{\g} a,b]_\g=[a,\dM_{\g} b]_\g$,
\item[\rm(c)] $[[x,y]_\g,z]_\g+[[y,z]_\g,x]_\g+[[z,x]_\g,y]_\g=0$,
\item[\rm(d)] $[[x,y]_\g,a]_\g+[[y,a]_\g,x]_\g+[[a,x]_\g,y]_\g=0$.
\end{itemize}
\end{defi}

We denote a strict Lie 2-algebra by $(\g,\dM_{\g},[\cdot,\cdot]_\g)$.

Let $ \mathbb V:V_1\stackrel{\partial}{\longrightarrow} V_0$ be a 2-term
complex of vector spaces, we can form a new 2-term complex of vector
spaces $\End(\mathbb V):\End^1(\mathbb
V)\stackrel{\delta}{\longrightarrow} \End^0_\partial(\mathbb V)$ by
defining $\delta(D)=(\partial\circ D,D\circ\partial)$ for any
$D\in\End^1(\mathbb V)$, where $\End^1(\mathbb V)=\Hom(V_0,V_1)$ and
$$\End^0_\partial(\mathbb V)=\{X=(X_0,X_1)\in \End(V_0,V_0)\oplus \End(V_1,V_1)|~X_0\circ \partial=\partial\circ X_1\}.$$
Define  $[\cdot,\cdot]_C:\End^i(\mathbb V)\wedge\End^j(\mathbb V)\longrightarrow \End^{i+j}(\mathbb
V),0\le i+j\le 1$ by setting:
\begin{eqnarray}
[(X_0,X_1),(Y_0,Y_1)]_C&=&([X_0,Y_0],[X_1,Y_1]),\\
\label{no}[(X_0,X_1),D]_C&=&X_1\circ D-D\circ X_0.
\end{eqnarray}

\begin{thm}\label{thm:End(V)} With the above notations,
$(\End(\mathbb V),\delta,[\cdot,\cdot]_C)$ is a strict Lie $2$-algebra.
\end{thm}

\begin{defi}\label{morphism}
Let $(\g,\dM_\g,[\cdot,\cdot]_\g)$ and $(\g',\dM_{\g'},[\cdot,\cdot]_{\g'})$ be strict Lie $2$-algebras. A strict
Lie $2$-algebra homomorphism $f:\g\rightarrow \g'$ consists of
  two linear maps $f_{0}:\g_{0}\rightarrow \g_{0}'$ and $f_{1}:\g_{1}\rightarrow \g_{1}',$
 such that the following equalities hold for all $ x,y\in \g_{0},
a\in \g_{1},$
\begin{itemize}
\item[$\bullet$] $\dM_{\g'}\circ f_1=f_0\circ\dM_\g$,
\item[$\bullet$] $f_{0}[x,y]_\g=[f_{0}(x),f_{0}(y)]_{\g'},$
\item[$\bullet$] $f_{1}[x,a]_\g=[f_{0}(x),f_{1}(a)]_{\g'}$.
\end{itemize}
\end{defi}

A strict representation of a strict Lie 2-algebra $(\g,\dM_\frkg,[\cdot,\cdot]_\frkg)$ on a 2-term complex $\mathbb V$ is a strict Lie 2-algebra homomorphism
$\rho=(\rho_0,\rho_1):\g\lon \End(\mathbb V)$. Given a  strict representation, we have the corresponding generalized Chevalley-Eilenberg complex
$(C^i(\g;\mathbb V),\D_\rho)$, where $C^i(\g;\mathbb V)$
denote the set of strict Lie 2-algebra $i$-cochains defined by\footnote{Here the degree of an element in $\Hom(\wedge^p\g_0\wedge\odot^q\g_1,V_s)$ where for $s=0,1$ is $p+2q-s$.}
$$C^i(\g;\mathbb V)=\bigoplus_{p+2q-s=i\atop p\neq i+1}\Hom(\wedge^p\g_0\wedge\odot^q\g_1,V_s),$$
and the coboundary operator $\D_{\rho}=\hat{\dM_\g}+\hat{\partial}+\dM_{\rho}^{(0,1)}+\dM_{\rho}^{(1,0)},$
in which
\begin{eqnarray*}
\hat{\dM_\g}:\Hom(\wedge^p\g_0\wedge\odot^q\g_1,V_s)&\lon&\Hom(\wedge^{p-1}\g_0\wedge\odot^{q+1}\g_1,V_s),\\
\hat{\partial}:\Hom(\wedge^p\g_0\wedge\odot^q\g_1,V_1)&\lon&\Hom(\wedge^{p}\g_0\wedge\odot^{q}\g_1,V_0),\\
\dM_{\rho}^{(0,1)}:\Hom(\wedge^p\g_0\wedge\odot^q\g_1,V_0)&\lon&\Hom(\wedge^{p}\g_0\wedge\odot^{q+1}\g_1,V_1),\\
\dM_{\rho}^{(1,0)}:\Hom(\wedge^p\g_0\wedge\odot^q\g_1,V_s)&\lon&\Hom(\wedge^{p+1}\g_0\wedge\odot^{q}\g_1,V_s),
\end{eqnarray*}
are defined by
\begin{eqnarray*}
(\hat{\dM_\g}f)(x_1,\cdots,x_{p-1},a_1,\cdots,a_{q+1})&=&(-1)^p\big(f(x_1,\cdots,x_{p-1},\dM_\g a_1,\cdots,a_{q+1})\\
                                                       &&+c.p.(a_1,\cdots,a_{q+1})\big),\\
                 \hat{\partial}f&=&(-1)^{p+2q}\partial\circ f, \\
(\dM_{\rho}^{(0,1)}f)(x_1,\cdots,x_{p},a_1,\cdots,a_{q+1})&=&(-1)^p\sum_{i=1}^{q+1}\rho_1(a_i)f(x_1,\cdots,x_{p},a_1,\cdots,\widehat{a_i},\cdots,a_{q+1}),
\end{eqnarray*}
and
\begin{eqnarray*}
  &&(\dM_{\rho}^{(1,0)}f)(x_1,\cdots,x_{p+1},a_1,\cdots,a_{q})=\sum_{i=1}^{p+1}(-1)^{i+1}\rho_0(x_i)f(x_1,\cdots,\widehat{x_i},\cdots,x_{p+1},a_1,\cdots,a_{q})\\
                                                          &&\qquad\qquad+\sum_{i<j}(-1)^{i+j}f([x_i,x_j]_\g,x_1,\cdots,\widehat{x_i},\cdots,\widehat{x_j},\cdots,x_{p+1},a_1,\cdots,a_{q})\\
                                                          &&\qquad\qquad+\sum_{i,j}(-1)^{i}f(x_1,\cdots,\widehat{x_i},\cdots,x_{p+1},a_1,\cdots,[x_i,a_j]_\g,\cdots,a_{q}),\\
\end{eqnarray*}
for all $x_i\in\g_0,a_i\in\g_1,i\in\mathbb N$. Denote the set of cocycles by $\huaZ(\g;\mathbb V)$ and the set of coboundaries by $\huaB(\g;\mathbb V)$. The corresponding cohomology is denoted by $\huaH(\g;\mathbb V)_{\rho}$.
\begin{rmk}\label{strict cohomology}
The cohomology of a strict Lie $2$-algebra with the coefficient in a strict representation $\mathbb V$ was given in \cite{BSZ}. See \cite{LM,LSZ} for more details for the general case. The cohomology complex given above is its subcomplex.
\end{rmk}

\emptycomment{
is given by
\begin{itemize}
\item[]degree 0: $V_0\stackrel{\D}{\longrightarrow}$
\item[]degree 1: $\Hom(\g_0,V_0)\oplus\Hom(\g_1,V_1)\stackrel{\D}{\longrightarrow}$
\item[]degree 2: $\Hom(\g_1,V_0)\oplus\Hom(\wedge^2\g_0,V_0)\oplus\Hom(\g_0\wedge\g_1,V_1)\stackrel{\D}{\longrightarrow}$
\item[]degree 3: $\Hom(\g_0\wedge\g_1,V_0)\oplus\Hom(\odot^2\g_1,V_1)\oplus\Hom(\wedge^3\g_0,V_0)\oplus\Hom(\wedge^2\g_0\wedge\g_1,V_1)\stackrel{\D}{\longrightarrow}\cdots$
\end{itemize}
}

\section{Strict derivations on strict Lie 2-algebras}

In this section, we outline some basic results about strict derivations, inner derivations, outer derivations and center of strict Lie 2-algebras. It is not totally new. See \cite{CSZ,DL,LL,LLS} for more details.

\begin{defi}\label{Def:derivation}
Let $(\g,\dM_\frkg,[\cdot,\cdot]_\frkg)$ be a strict Lie $2$-algebra. A strict derivation of degree $0$
of $\frkg$ is an element $(X_0,X_1)\in \End_{\dM_\g}^{0}(\g)$, such that for all $x,y\in\g_0$ and $a\in\g_1$,
\begin{eqnarray}
\label{der1}X_0[x,y]_\frkg&=&[X_0x,y]_\frkg+[x,X_0y]_\frkg,\\
\label{der2}X_1[x,a]_\frkg&=&[X_0x,a]_\frkg+[x,X_1a]_\frkg.
\end{eqnarray}
\end{defi}
We denote the set of strict derivations of degree 0 of $\g$ by $\SDer^0(\frkg)$.

\begin{defi}\label{derivation1}
Let $(\g,\dM_\frkg,[\cdot,\cdot]_\frkg)$ be a strict Lie $2$-algebra. A strict derivation of degree $1$
of $\frkg$ is a linear map $\Theta:\g_0\lon\g_1$ such that for all $x,y\in\g_0$,
\begin{eqnarray}
\label{der3}\Theta[x,y]_\g=[\Theta(x),y]_\g+[x,\Theta(y)]_\g.
\end{eqnarray}
\end{defi}
We denote the set of strict derivations of degree 1 of $\g$ by $\SDer^1(\frkg)$.

It is straightforward to obtain that

\begin{lem}\label{2term}
For all $\Theta\in\SDer^1(\frkg)$, we have $\delta(\Theta)\in \SDer^0(\frkg).$ Thus, we have a well-defined complex:
\begin{eqnarray}
\SDer(\frkg):\SDer^1(\frkg)\stackrel{\delta}{\longrightarrow} \SDer^0(\frkg).
\end{eqnarray}
\end{lem}

 By straightforward computations, we have
\begin{lem}\label{bracket}
For all $(X_0,X_1),(Y_0,Y_1)\in\SDer^0(\frkg)$, $\Theta\in\SDer^1(\frkg)$, we have
\begin{eqnarray*}
[(X_0,X_1),(Y_0,Y_1)]_C\in\SDer^0(\frkg),\,\,\,\,\,[(X_0,X_1),\Theta]_C\in\SDer^1(\frkg).
\end{eqnarray*}
\end{lem}

By lemma \ref{2term} and \ref{bracket}, we have
\begin{thm}\label{derivation Lie 2-algebra}
With the above notations, $(\SDer(\frkg),\delta,[\cdot,\cdot]_C)$ is a strict Lie $2$-algebra, which is a sub-algebra of the strict Lie $2$-algebra $(\End(\g),\delta,[\cdot,\cdot]_C)$. We call it the {\bf derivation Lie $2$-algebra} of the strict Lie $2$-algebra $(\g,\dM_\frkg,[\cdot,\cdot]_\frkg)$.
\end{thm}

\begin{rmk}
Our derivation Lie $2$-algebra $\SDer(\frkg)$ is equivalent to the actor of a Lie algebra crossed module  given in \cite{CL}.
\end{rmk}

For any strict Lie 2-algebra $\g$, there is a strict Lie 2-algebra homomorphism from $\g$ to $\SDer(\frkg)$ given as follows. Define
$\ad_0:\g_0\lon \SDer^0(\frkg)$ and $\ad_1:\g_1\lon \SDer^1(\frkg)$
by
$$\ad_0(x)=(\ad^0_x,\ad^1_x),\,\,\,\,\ad_1(a)(x)=[a,x]_\g,$$
where $\ad^0:\g_0\longrightarrow \gl(\g_0)$ and $\ad^1:\g_0\longrightarrow \gl(\g_1)$ are defined by
$$
\ad^0_xy=[x,y]_\g,\quad \ad^1_xa=[x,a]_\g.
$$

\begin{lem}\label{adjoint}
With the above notations, $\frkad=(\ad_0,\ad_1)$ is a strict Lie $2$-algebra homomorphism from $\g$ to $\SDer(\frkg)$, which is called the {\bf adjoint representation} of $\g$.
\end{lem}

\pf It is straightforward. \qed

\begin{defi}\label{def:center}
Let $(\g,\dM_\frkg,[\cdot,\cdot]_\frkg)$ be a strict Lie $2$-algebra. The center of $(\g,\dM_\frkg,[\cdot,\cdot]_\frkg)$, denoted by $\cen(\g)$, is defined as the kernel  of the strict Lie $2$-algebra homomorphism $\frkad:\g\lon \SDer(\g).$
\end{defi}

We write $\cen(\g)=\cen^1(\g)\oplus \cen^0(\g)$, where $\cen^1(\g)=\ker(\ad_1)$ and $\cen^0(\g)=\ker(\ad_0).$
 It is obvious that
\begin{lem}\label{lem:center}
The $\cen(\g)$ is an ideal of the strict Lie $2$-algebra $(\g,\dM_\frkg,[\cdot,\cdot]_\frkg)$.
\end{lem}

\begin{pro}\label{thm:center}
Let $(\g,\dM_\frkg,[\cdot,\cdot]_\frkg)$ be a strict Lie $2$-algebra. Then we have
$$\huaH^0(\g;\g)_{\frkad}=\cen^0(\g).$$
\end{pro}

\pf For all $x\in\g_0$, $\D_{\frkad}(x)$ has two components as follows:
\begin{equation}
\left\{
\label{cob1}\begin{array}{lcll}
{\Lambda_1} &=&\dM_{\rho}^{(1,0)}(x) &\in\Hom(\mathfrak{g}_{0},\mathfrak{g}_{0}),\\
{\Lambda_2} &=&\dM_{\rho}^{(0,1)}(x) &\in\Hom(\g_1,\g_1).
\end{array}
\right.
\end {equation}
For all $y\in\g_0,b\in\g_1$, we have
$$
 \Lambda_1(y)=[y,x]_\g,\quad
 \Lambda_2(b)=[b,x]_\g.
$$
Thus,  $\D_{\frkad}(x)=0$ if and only if $x\in\cen^0(\g)$.   \qed\vspace{3mm}

Since $\frkad$ is a strict Lie $2$-algebra homomorphism, it is obvious that $\Img(\frkad$) is a sub-algebra of the strict Lie $2$-algebra $\SDer(\frkg)$. Then we can get a strict Lie $2$-algebra $\SInn(\frkg)$ given by
\begin{eqnarray}
\SInn(\frkg):\SInn^1(\frkg)\stackrel{\triangle}=\Img(\ad_1)\stackrel{\delta}{\longrightarrow} \SInn^0(\frkg)\stackrel{\triangle}=\Img(\ad_0).
\end{eqnarray}

\begin{lem}\label{ideal}
Let $(\g,\dM_\frkg,[\cdot,\cdot]_\frkg)$ be a strict Lie $2$-algebra. For all $x\in\g_0,a\in\g_1$ and $(X_0,X_1)\in\SDer^0(\frkg),\Theta\in\SDer^1(\frkg)$, we have
$$[(X_0,X_1),\ad_1(a)]_C=\ad_1(X_1a),\,\,[(X_0,X_1),\ad_0(x)]_C=\ad_0(X_0x),\,\,[\Theta,\ad_0(x)]_C=\ad_1(\Theta x).$$
Therefore, $\SInn(\frkg)$ is an ideal of the strict Lie $2$-algebra $\SDer(\frkg)$.
\end{lem}

\pf It follows by straightforward computations. \qed\vspace{3mm}

Denote by $\SOut(\g)$ the set of  out strict derivations of the strict Lie $2$-algebra $(\g,\dM_\frkg,[\cdot,\cdot]_\frkg)$, i.e.
$$\SOut(\g)=\SDer(\g)/\SInn(\g).$$
We use $\pi=(\pi_0,\pi_1)$ to denote the quotient map from $\SDer(\h)$ to $\SOut(\h)$.

\begin{pro}\label{outer}
Let $(\g,\dM_\frkg,[\cdot,\cdot]_\frkg)$ be a strict Lie $2$-algebra. We have
$$\huaH^1(\g;\g)_{\frkad}=\SOut^0(\g).$$
\end{pro}

\pf For all $X_0\in\Hom(\g_0,\g_0)$ and $X_1\in\Hom(\g_1,\g_1)$,  $\D_{\frkad}(X_0,X_1)$ has three components as follows:

\begin{equation}
\left\{
\begin{array}{lcll}
{\Gamma_1} &=&\hat{\dM_\g}X_0+\hat{\partial}X_1                 &\in\Hom(\mathfrak{g}_{1},\mathfrak{g}_{0}),\\
{\Gamma_2} &=&\dM_{\rho}^{(1,0)}X_0                             &\in\Hom(\wedge^2\g_0,\g_0),\\
{\Gamma_2} &=&\dM_{\rho}^{(0,1)}X_0+\dM_{\rho}^{(1,0)}X_1       &\in\Hom(\g_0\wedge\g_1,\g_1).
\end{array}
\right.
\end{equation}
For $x,y\in\g_0,a\in\g_1$, we have
\begin{equation*}
\left\{
\begin{array}{lcl}
{\Gamma_1}(a) &=&\dM_{\g}(X_1a)-X_0(\dM_\g a),\\
{\Gamma_2}(x,y) &=&[x,X_0y]_\g-[y,X_0x]_\g-X_0[x,y]_\g,\\
{\Gamma_3}(x,a) &=&[x,X_1a]_\g-X_1[x,a]_\g-[a,X_0x]_\g.\\
\end{array}
\right.
\end{equation*}
Thus,   $(X_0,X_1)\in \huaZ^1(\g;\g)$ if and only if $(X_0,X_1)\in\SDer^0(\g)$. Furthermore, it is straightforward to deduce that  $(X_0,X_1)\in \huaB^1(\g;\g)$ if and only if $(X_0,X_1)\in\SInn^0(\g)$. The proof is finished. \qed

\section{Non-abelian extensions of strict Lie 2-algebras}

In this section, we give a general description of a non-abelian extension by choosing a section.
 \begin{defi}\label{eq:ext}
 \begin{itemize}
 \item[\rm (i)] Let $(\g,\dM_{\g},[\cdot,\cdot]_\g)$, $(\h,\dM_{\h},[\cdot,\cdot]_\h)$, $(\hat{\mathfrak{g}},\hat{\dM},[\cdot,\cdot]_{\hat{\g}})$ be Lie $2$-algebras and
$i=(i_{0},i_{1}):\frkh\longrightarrow\hat{\mathfrak{g}},~~p=(p_{0},p_{1}):\hat{\mathfrak{g}}\longrightarrow\mathfrak{g}$
be strict homomorphisms. The following sequence of Lie $2$-algebras is a
short exact sequence if $\mathrm{Im}(i)=\ker(p)$,
$\ker(i)=0$ and $\mathrm{Im}(p)=\g$,
\begin{equation}\label{eq:ext1}
\CD
  0     @>>>            \frkh_1 @>i_1>>                 \hat{\g_1} @>p_1>>          \g_1 @>>>     0 \\
  @.                       @V \dM_\frkh VV                  @V \hat{\dM} VV           @V\dM_\g VV   @.      \\
  0     @>>>            \frkh_{0} @>i_0>>               \hat{\g_0} @>p_0>>          \g_0@>>>      0.
\endCD
\end{equation}

We call $\hat{\mathfrak{g}}$  a non-abelian extension of $\mathfrak{g}$ by
$\frkh$, and denote it by $\E_{\hat{\g}}.$
\item[\rm (ii)] A section $\sigma:\mathfrak{g}\longrightarrow\hat{\mathfrak{g}}$ of $p:\hat{\mathfrak{g}}\longrightarrow\mathfrak{g}$
consists of linear maps $\sigma_0:\mathfrak{g}_0\longrightarrow\hat{\g_0}$ and $\sigma_1:\mathfrak{g}_1\longrightarrow\hat{\g_1}$
 such that  $p_0\circ\sigma_0=id_{\mathfrak{g}_0}$ and  $p_1\circ\sigma_1=id_{\mathfrak{g}_1}$.
\item[\rm (iii)] We say two non-abelian extensions of Lie $2$-algebras
 $\E_{\hat{\g}}:\frkh\stackrel{i}{\longrightarrow}\hat{\g}\stackrel{p}{\longrightarrow}\g$
 and $\E_{\tilde{\g}}:\frkh\stackrel{j}{\longrightarrow}\tilde{\g}\stackrel{q}{\longrightarrow}\g$ are isomorphic
 if there exists a strict Lie $2$-algebra homomorphism $F:\hat{\g}\longrightarrow\tilde{\g}$  such that $F\circ i=j$, $q\circ
 F=p$.
\end{itemize}
\end{defi}

It is easy to see that this implies that $F$ is an isomorphism of strict Lie $2$-algebras, hence defines an equivalence relation.

  Since $(i_0,i_1)$ are inclusions and $(p_0,p_1)$ are projections, a section $\sigma$ induces linear maps:
$$
\begin{array}{rlclcrcl}
\varphi:&\mathfrak{g}_{1}&\longrightarrow&\frkh_{0},&& \varphi(a)&=&\hat{\dM}\sigma_1(a)-\sigma_0(\dM_\g a),\\
\mu^{0}:&\mathfrak{g}_{0}&\longrightarrow& \gl(\frkh_0),&& \mu^{0}(x)(u)&=&[\sigma_0(x),u]_{\hat{\mathfrak{g}}},\\
\mu^{1}:&\mathfrak{g}_{0}&\longrightarrow& \gl(\frkh_1),&& \mu^{1}(x)(m)&=&[\sigma_0(x),m]_{\hat{\mathfrak{g}}},\\
\mu_{1}:&\mathfrak{g}_{1}&\longrightarrow& \Hom(\frkh_0,\frkh_1),&&\mu_{1}(a)(u)&=&[\sigma_1(a),u]_{\hat{\mathfrak{g}}},\\
\omega:&\wedge^2\mathfrak{g}_{0}&\longrightarrow&\frkh_{0},&& \omega(x,y)&=&[\sigma_0(x),\sigma_0(y)]_{\hat{\mathfrak{g}}}-\sigma_0[x,y]_\mathfrak{g},\\
\nu:&\mathfrak{g}_{0}\wedge\mathfrak{g}_{1}&\longrightarrow&\frkh_{1},&& \nu(x,a)&=&[\sigma_0(x),\sigma_1(a)]_{\hat{\mathfrak{g}}}-\sigma_1[x,a]_\mathfrak{g},
\end{array}
$$
for all $x,y,z\in\mathfrak{g}_{0}$, $a,b\in\mathfrak{g}_{1}$,
$u,v\in\frkh_{0}$ and $m,n\in\frkh_{1}$.

Obviously, $\hat{\frkg}$ is isomorphic to $\frkg\oplus\frkh$ as vector spaces. Transfer the strict Lie 2-algebra structure on $\hat{\frkg}$
to that on $\frkg\oplus\frkh$, we obtain a strict Lie 2-algebra $(\frkg\oplus\frkh, \dM_{\frkg\oplus\frkh}, [\cdot, \cdot]_{\frkg\oplus\frkh})$, where $\dM_{\frkg\oplus\frkh}$ and $[\cdot, \cdot]_{\frkg\oplus\frkh}$ are given by
\begin{eqnarray}
\label{d}\dM_{\frkg\oplus\frkh}(a+m) &=&\dM_{\frkg}a+\varphi(a)+\dM_{\frkh}(m),\\
\label{bo}[x+u,y+v]_{\frkg\oplus\frkh}&=&[x,y]_\g+\omega(x,y)+\mu^{0}(x)v-\mu^{0}(y)u+[u,v]_{\frkh},\\
\label{b1}[x+u,a+m]_{\frkg\oplus\frkh}&=&[x,a]_\g+\nu(x,a)+\mu^{1}(x)m-\mu_{1}(a)u+[u,m]_{\frkh}.
\end{eqnarray}
The following proposition gives the conditions on $\varphi,\mu^{0},\mu^{1},\mu_{1},\omega$ and $\nu$ such that $(\frkg\oplus\frkh, \dM_{\frkg\oplus\frkh}, [\cdot, \cdot]_{\frkg\oplus\frkh})$ is a strict Lie 2-algebra. For convenience, we denote by $\mu_0=(\mu^0,\mu^1)$.

\begin{pro}\label{extension}
With the above notations, $(\frkg\oplus\frkh, \dM_{\frkg\oplus\frkh}, [\cdot, \cdot]_{\frkg\oplus\frkh})$ a strict Lie $2$-algebra if and only
if the following equalities hold:
\begin{eqnarray}
\emptycomment{
\label{p1}\dM_\frkh\circ\mu^{1}(x)                &=&\mu^{0}(x)\circ\dM_\frkh,\\
\label{p2}\mu^{0}(x)[u,v]_\frkh                   &=&[\mu^{0}(x)u,v]_\frkh+[u,\mu^{0}(x)v]_\frkh,\\
\label{p3}\mu^{1}(x)[m,u]_\frkh                  &=&[\mu^{1}(x)m,u]_\frkh+[m,\mu^{0}(x)u]_\frkh\\
\label{p4}\mu_{1}(a)[u,v]_\frkh                  &=&[\mu_{1}(a)u,v]_\frkh+[u,\mu_{1}(a)v]_\frkh,\\
\label{p5}\delta(\mu_1(a))                                  &=&\mu_0(\dM_\frkg(a))+\ad_0(\varphi(a)),\\
\label{p6}[\mu_0(x),\mu_0(y)]_C                    &=&\mu_0([x,y]_\g)+\ad_0(\omega(x,y)),\\
\label{p7}[\mu_1(a),\mu_0(x)]_C                    &=&\mu_1([a,x]_\g)+\ad_1(\nu(a,x)),\\
\label{p8}\varphi([x,a]_\frkg)+\dM_\frkh(\nu(x,a))&=&\omega(x,\dM_\frkg(a))+\mu^{0}(x)\varphi(a),\\
\label{p9}\nu(\dM_\frkg(a),b)-\mu_{1}(b)\varphi(a)&=&\nu(a,\dM_\frkg(b))+\mu_{1}(a)\varphi(b),\\
\label{p10}\mu^{0}(x)\omega(y,z)+c.p.&=&\omega([x,y]_\frkg,z)+c.p.,\\
\label{p11}\nu([x,y]_\frkg,a)+c.p.&=&\mu^{1}(x)\nu(y,a)+\mu^{1}(y)\nu(a,x)+\mu_{1}(a)\omega(x,y).\\
}
\label{p1}\dM_\frkh\circ\mu^{1}(x)                &=&\mu^{0}(x)\circ\dM_\frkh,\\
\label{p2}\mu^{0}(x)[u,v]_\frkh                   &=&[\mu^{0}(x)u,v]_\frkh+[u,\mu^{0}(x)v]_\frkh,\\
\label{p3}\mu^{1}(x)[m,u]_\frkh                  &=&[\mu^{1}(x)m,u]_\frkh+[m,\mu^{0}(x)u]_\frkh,\\
\label{p4}\mu_{1}(a)[u,v]_\frkh                  &=&[\mu_{1}(a)u,v]_\frkh+[u,\mu_{1}(a)v]_\frkh,\\
\label{p5}\dM_\frkh\circ\mu_{1}(a)                &=&\mu^{0}(\dM_\frkg(a))+\ad^0_{\varphi(a)},\\
\label{p6}\mu_{1}(a)\circ \dM_\frkh               &=&\mu^{1}(\dM_\frkg(a))+\ad^1_{\varphi(a)},\\
\label{p7}[\mu^{0}(x),\mu^{0}(y)]&=&\mu^{0}([x,y]_\frkg)+\ad^0_{\omega(x,y)},\\
\label{p8}[\mu^{1}(x),\mu^{1}(y)]&=&\mu^{1}([x,y]_\frkg)+\ad^1_{\omega(x,y)},\\
\label{p9}\mu_{1}(a)\circ\mu^{0}(x)-\mu^{1}(x)\circ\mu_{1}(a)&=&\mu_{1}([a,x]_\frkg)+\ad_1(\nu(a,x)),\\
\label{p10}\varphi([x,a]_\frkg)+\dM_\frkh(\nu(x,a))&=&\omega(x,\dM_\frkg a)+\mu^{0}(x)\varphi(a),\\
\label{p11}\nu(\dM_\frkg a,b)-\mu_{1}(b)\varphi(a)&=&\nu(a,\dM_\frkg b)+\mu_{1}(a)\varphi(b),\\
\label{p12}\mu^{0}(x)\omega(y,z)+c.p.&=&\omega([x,y]_\frkg,z)+c.p.,\\
\label{p13}\nu([x,y]_\frkg,a)+c.p.&=&\mu^{1}(x)\nu(y,a)+\mu^{1}(y)\nu(a,x)+\mu_{1}(a)\omega(x,y).
\end{eqnarray}
\end{pro}

\pf Let $(\frkg\oplus\frkh, \dM_{\frkg\oplus\frkh}, [\cdot, \cdot]_{\frkg\oplus\frkh})$ be a strict Lie 2-algebra. By
$$\dM_{\frkg\oplus\frkh}[x+u,a+m]_{\frkg\oplus\frkh}=[x+u,\dM_{\frkg\oplus\frkh}(a+m)]_{\frkg\oplus\frkh},$$
we deduce that \eqref{p1}, \eqref{p5} and \eqref{p10} hold. By
$$[\dM_{\frkg\oplus\frkh}(a+m),b+n]_{\frkg\oplus\frkh}=[a+m,\dM_{\frkg\oplus\frkh}(b+n)]_{\frkg\oplus\frkh},$$
we deduce that \eqref{p6} and \eqref{p11} hold. By
$$[[u,v]_{\frkg\oplus\frkh},x]_{\frkg\oplus\frkh}+[[v,x]_{\frkg\oplus\frkh},u]_{\frkg\oplus\frkh}+[[x,u]_{\frkg\oplus\frkh},v]_{\frkg\oplus\frkh}=0,$$
we deduce that \eqref{p2} holds.  By
$$[[u,x]_{\frkg\oplus\frkh},m]_{\frkg\oplus\frkh}+[[x,m]_{\frkg\oplus\frkh},u]_{\frkg\oplus\frkh}+[[m,u]_{\frkg\oplus\frkh},x]_{\frkg\oplus\frkh}=0,$$
we deduce that \eqref{p3} holds. By
$$[[u,v]_{\frkg\oplus\frkh},a]_{\frkg\oplus\frkh}+[[v,a]_{\frkg\oplus\frkh},u]_{\frkg\oplus\frkh}+[[a,u]_{\frkg\oplus\frkh},v]_{\frkg\oplus\frkh}=0,$$
we deduce that \eqref{p4} holds. By
$$[[u,x]_{\frkg\oplus\frkh},y]_{\frkg\oplus\frkh}+[[x,y]_{\frkg\oplus\frkh},u]_{\frkg\oplus\frkh}+[[y,u]_{\frkg\oplus\frkh},x]_{\frkg\oplus\frkh}=0,$$
we deduce that \eqref{p7} holds. By
$$[[x,y]_{\frkg\oplus\frkh},m]_{\frkg\oplus\frkh}+[[y,m]_{\frkg\oplus\frkh},x]_{\frkg\oplus\frkh}+[[m,x]_{\frkg\oplus\frkh},y]_{\frkg\oplus\frkh}=0,$$
we deduce that \eqref{p8} holds. By
$$[[u,x]_{\frkg\oplus\frkh},a]_{\frkg\oplus\frkh}+[[x,a]_{\frkg\oplus\frkh},u]_{\frkg\oplus\frkh}+[[a,u]_{\frkg\oplus\frkh},x]_{\frkg\oplus\frkh}=0,$$
we deduce that \eqref{p9} holds.  By
$$[[x,y]_{\frkg\oplus\frkh},z]_{\frkg\oplus\frkh}+[[y,z]_{\frkg\oplus\frkh},x]_{\frkg\oplus\frkh}+[[z,x]_{\frkg\oplus\frkh},y]_{\frkg\oplus\frkh}=0,$$
we deduce that \eqref{p12} holds. By
$$[[x,y]_{\frkg\oplus\frkh},a]_{\frkg\oplus\frkh}+[[y,a]_{\frkg\oplus\frkh},x]_{\frkg\oplus\frkh}+[[a,x]_{\frkg\oplus\frkh},y]_{\frkg\oplus\frkh}=0,$$
we deduce that \eqref{p13} holds.

Conversely, if \eqref{p1}-\eqref{p13} holds, it is straightforward to see that $(\frkg\oplus\frkh, \dM_{\frkg\oplus\frkh}, [\cdot, \cdot]_{\frkg\oplus\frkh})$ is a strict Lie 2-algebra. The proof is finished. \qed\vspace{3mm}

For any non-abelian extension, by choosing a section, it is isomorphic to $(\frkg\oplus\frkh, \dM_{\frkg\oplus\frkh}, [\cdot, \cdot]_{\frkg\oplus\frkh})$. Therefore, we only consider non-abelian extensions of the form $(\frkg\oplus\frkh, \dM_{\frkg\oplus\frkh}, [\cdot, \cdot]_{\frkg\oplus\frkh})$ in the sequel.

\section{Classification of non-abelian extensions of strict Lie 2-algebras: special case}

In this section, we classify non-abelian extensions of strict Lie 2-algebras for the case that $\cen(\h)=0.$

\begin{pro}\label{pro:special}
Let $(\g,\dM_\frkg,[\cdot,\cdot]_\frkg)$ and $(\frkh,\dM_\frkh,[\cdot,\cdot]_\frkh)$ be strict Lie $2$-algebras such that $\cen(\frkh)=0$. Then isomorphism classes of non-abelian extensions of $\g$ by $\frkh$ correspond bijectively to strict Lie $2$-algebra homomorphisms
\begin{eqnarray*}
\bar{\mu}:\g\lon \SOut(\frkh).
\end{eqnarray*}
\end{pro}

\pf Let $(\frkg\oplus\frkh, \dM_{\frkg\oplus\frkh}, [\cdot, \cdot]_{\frkg\oplus\frkh})$ be a non-abelian extension of $\g$ by $\frkh$ given by \eqref{d}-\eqref{b1}. For $x\in\g_0,a\in\g_1$, we denote by $\mu_0(x)=(\mu^{0}(x),\mu^{1}(x)).$
By \eqref{p1}-\eqref{p3}, we have $\mu_0(x)\in\SDer^0(\frkh)$. By
\eqref{p4}, we have $\mu_1(a)\in\SDer^1(\frkh)$. Let  $\pi=(\pi_1,\pi_0):\SDer(\frkh)\lon \SOut(\frkh)$ be  the quotient map. We denote  the induced strict Lie 2-algebra structure on $\SOut(\frkh)$ by $[\cdot,\cdot]_C'$ and $\delta'$. Hence we can define $$\bar{\mu}=(\bar{\mu}_0,\bar{\mu}_1),\,\quad\mbox{where}\quad\,\bar{\mu}_0=\pi_0\circ\mu_0,\quad \bar{\mu}_1=\pi_1\circ\mu_1.$$
 By \eqref{p5} and \eqref{p6}, we have
\begin{eqnarray*}
\delta'(\bar{\mu}_1(a))=\delta'(\pi_1(\mu_{1}(a)))=\pi_0(\delta(\mu_{1}(a)))&=&\pi_0(\dM_\frkh\circ\mu_{1}(a),\mu_{1}(a)\circ\dM_\frkh)\\
                                                       &=&\pi_0(\mu^0(\dM_\frkg a)+\ad^0_{\varphi(a)},\mu^1(\dM_\frkg a)+\ad^1_{\varphi(a)})\\
                                                       &=&\pi_0(\mu_0(\dM_\frkg a)+\ad_0(\varphi(a)))\\
                                                       &=&\bar{\mu}_0(\dM_\frkg a).
\end{eqnarray*}
By \eqref{p7} and \eqref{p8}, we have
\begin{eqnarray*}
\bar{\mu}_0([x,y]_\g)=(\pi_0\circ\mu_0)([x,y]_\g)&=&\pi_0(\mu^0([x,y]_\g),\mu^1([x,y]_\g))\\
                      &=&\pi_0([\mu^0(x),\mu^0(y)]-\ad^0_{\omega(x,y)},[\mu^1(x),\mu^1(y)]-\ad^1_{\omega(x,y)})\\
                      &=&\pi_0([\mu_0(x),\mu_0(y)]_C-\ad_0(\omega(x,y)))\\
                       &=&[\bar{\mu}_0(x),\bar{\mu}_0(y)]_C'.
\end{eqnarray*}
Similarly, by \eqref{p9}, we have
\begin{eqnarray*}
\bar{\mu}_1([a,x]_\g)                      &=&[\bar{\mu}_1(a),\bar{\mu}_0(x)]_C'.
\end{eqnarray*}
Thus, $\bar{\mu}$ is a strict Lie $2$-algebra homomorphism from $\g$ to $\SOut(\frkh)$.

Let $(\frkg\oplus\frkh, \dM_{\frkg\oplus\frkh}, [\cdot, \cdot]_{\frkg\oplus\frkh})$ and $(\frkg\oplus\frkh, \dM_{\frkg\oplus\frkh}', [\cdot, \cdot]_{\frkg\oplus\frkh}')$ be isomorphic extensions of $\g$ by $\frkh$. Then there is a strict Lie $2$-algebra homomorphism $\theta=(\theta_0,\theta_1):(\frkg\oplus\frkh, \dM_{\frkg\oplus\frkh}', [\cdot, \cdot]_{\frkg\oplus\frkh}')\lon(\frkg\oplus\frkh, \dM_{\frkg\oplus\frkh}, [\cdot, \cdot]_{\frkg\oplus\frkh})$, such that we have the following commutative diagram:
\[\begin{CD}
0@>>>(\frkh,\dM_\frkh,[\cdot,\cdot]_\frkh)@>\iota >>(\frkg\oplus\frkh, \dM_{\frkg\oplus\frkh}', [\cdot, \cdot]_{\frkg\oplus\frkh}')@>\pr>>(\g,\dM_\frkg,[\cdot,\cdot]_\frkg)  @>>>0\\
@.    @|                       @V\theta VV                                                   @|               @.\\
0@>>>(\frkh,\dM_\frkh,[\cdot,\cdot]_\frkh)@>\iota >>(\frkg\oplus\frkh, \dM_{\frkg\oplus\frkh}, [\cdot, \cdot]_{\frkg\oplus\frkh})@>\pr>>(\g,\dM_\frkg,[\cdot,\cdot]_\frkg)  @>>>0,
\end{CD}\]
where $\iota$ is the inclusion and $\pr$ is the projection. Since for all $x\in\g_0$ and $a\in\g_1$, we have $\pr_0(\theta_0(x))=x$ and $\pr_1(\theta_1(a))=a$, we can assume that
 $$\theta_0(x+u)=x+\xi(x)+u,\quad \theta_1(a+m)=a+\eta(a)+m,$$
 for some linear maps $\xi:\g_0\lon\frkh_0$ and $\eta:\g_1\lon\frkh_1$.

By the definition of  homomorphisms between strict Lie $2$-algebras, we get
\begin{eqnarray}
\label{iso1}\mu'_{0}(x)-\mu_0(x)&=&\ad_0({\xi(x)}),\\
\label{iso2}\mu'_1(a)-\mu_1(a)&=&\ad_1({\eta(a)}),\\
\label{iso3}\varphi'(a)-\varphi(a)&=&\dM_{\frkh}(\eta(a))-\xi(\dM_\g a),\\
\label{iso4}\omega'(x,y)-\omega(x,y)&=&\mu^0(x)(\xi(y))-\mu^0(y)(\xi(x))-\xi([x,y]_\g)+[\xi(x),\xi(y)]_{\frkh},\\
\label{iso5}\nu'(a,x)-\nu(a,x)&=&\mu_1(a)(\xi(x))-\mu^1(x)(\eta(a))-\eta([a,x]_\g)+[\eta(a),\xi(x)]_{\frkh}.
\end{eqnarray}
Therefore, by \eqref{iso1} we have
\begin{eqnarray*}
\bar{\mu}_0'(x)=\pi_0(\mu_0'(x))=\pi_0(\mu_0(x)+\ad_0({\xi(x)}))=\pi_0(\mu_0(x)) =\bar{\mu}_0(x).
\end{eqnarray*}
By \eqref{iso2}, we have
\begin{eqnarray*}
\bar{\mu}_1'(a)=\pi_1(\mu_1'(a))=\pi_1(\mu_1(a)+\ad_1({\eta(a)}))=\pi_1(\mu_1(a))=\bar{\mu}_1(a).
\end{eqnarray*}
Thus, we obtain that isomorphic non-abelian extensions of $(\g,\dM_\frkg,[\cdot,\cdot]_\frkg)$ by $(\frkh,\dM_\frkh,[\cdot,\cdot]_\frkh)$ correspond to the same strict Lie $2$-algebra homomorphism from $\g$ to $\SOut(\frkh)$.

Conversely, let $\bar{\mu}$ be a strict Lie $2$-algebra homomorphism from $\g$ to $\SOut(\frkh)$. By choosing a section $s=(s_0,s_1)$ of $\pi:\SDer(\g)\lon\SOut(\frkh)$, we define
\begin{eqnarray}
\label{lift1}\mu_0(x)&=&(s_0\circ\bar{\mu}_0)(x),\\
\label{lift2}\mu_1(a)&=&(s_1\circ\bar{\mu}_1)(a).
\end{eqnarray}
We have $\mu_0(x)=(\mu^{0}(x),\mu^{1}(x))\in\SDer^0(\frkh)$ and $\mu_1(a)\in\SDer^1(\frkh)$. Thus we get \eqref{p1}-\eqref{p4}. Since $\pi$ and $\bar{\mu}$ are strict Lie $2$-algebra homomorphisms, for all $a\in\g_1$, we have
\begin{eqnarray*}
\pi_0(\delta(\mu_1(a))-\mu_0(\dM_\g a))=\delta'(\pi_1(\mu_1(a)))-\pi_0(\mu_0(\dM_\g a))=\delta'(\bar{\mu}_1(a))-\bar{\mu}_0(\dM_\g a)=0.
\end{eqnarray*}
Thus, we obtain $\delta(\mu_1(a))-\mu_0(\dM_\g a)\in\SInn^0(\frkh)$. Since $\cen(\frkh)=0$, there is a unique linear map $\varphi:\g_1\lon\frkh_0$ such that
\begin{eqnarray}
\label{lift3}\delta(\mu_1(a))-\mu_0(\dM_\g a)=\ad_0({\varphi(a)}).
\end{eqnarray}
\emptycomment{
Moreover, we have
\begin{eqnarray*}
\delta(\mu_{2}(a))-(\mu_0(\dM_\g a),\mu_1(\dM_\g a))&=&(\dM_{\frkh}\circ\mu_{2}(a)-\mu_0(\dM_\g a),\mu_{2}(a)\circ\dM_{\frkh}-\mu_1(\dM_\g a))\\
                                                    &=&(\ad_{\varphi(a)},\ad_{\varphi(a)}).
\end{eqnarray*}
}
Thus, we obtain \eqref{p5} and \eqref{p6}. For all $x,y\in\g_0$, we have
\begin{eqnarray*}
\pi_0([\mu_0(x),\mu_0(y)]_C-\mu_0([x,y]_\g))&=&[\pi_0(\mu_0(x)),\pi_0(\mu_0(y))]_C'-\pi_0(\mu_0([x,y]_\g))\\
                                               &=&[\bar{\mu}_0(x),\bar{\mu}_0(y)]_C'-\bar{\mu}_0([x,y]_\g)\\
                                               &=&0,
\end{eqnarray*}
which implies  $[\mu_0(x),\mu_0(y)]_C-\mu_0([x,y]_\g)\in\SInn^0(\frkh)$.  Since $\cen(\frkh)=0$, there is a unique linear map
$\omega:\wedge^2\mathfrak{g}_{0}\lon\frkh_0$ such that
\begin{eqnarray}
\label{lift4}[\mu_0(x),\mu_0(y)]_C-\mu_0([x,y]_\g)=\ad_0({\omega(x,y)}).
\end{eqnarray}
\emptycomment{
Moreover, we have
\begin{eqnarray*}
&&[(\mu_0(x),\mu_1(x)),(\mu_0(y),\mu_1(y))]_C-(\mu_0([x,y]_\g),\mu_1([x,y]_\g))\\
&=&([\mu_0(x),\mu_0(y)]-\mu_0([x,y]_\g),[\mu_1(x),\mu_1(y)]-\mu_1([x,y]_\g))\\
&=&(\ad_{\omega(x,y)},\ad_{\omega(x,y)}).
\end{eqnarray*}
}
Thus, we obtain \eqref{p7} and \eqref{p8}. Similarly, there is a unique linear map
$\nu:\mathfrak{g}_{0}\wedge\mathfrak{g}_{1}\lon\frkh_1$ such that
\begin{eqnarray}
\label{lift5}[\mu_1(a),\mu_0(x)]_C-\mu_1([a,x]_\g)=\ad_1({\nu(a,x)}),
\end{eqnarray}
which implies that \eqref{p9} holds. For all $x\in\g_0,a\in\g_1,u\in\frkh_0$, by \eqref{lift3}-\eqref{lift5}, we have
\begin{eqnarray*}
&&[\varphi([x,a]_\frkg)+\dM_\frkh(\nu(x,a))-\omega(x,\dM_\frkg a)-\mu^{0}(x)\varphi(a),u]_{\frkh}\\
&=&\ad^0_{\varphi([x,a]_\frkg)}(u)-\dM_{\frkh}[\nu(a,x),u]_{\frkh}-\ad^0_{\omega(x,\dM_\frkg a)}(u)-(\mu^{0}(x)[\varphi(a),u]_{\frkh}-[\varphi(a),\mu^{0}(x)u])\\
&=&\big(\dM_\frkh\circ\mu_{1}([x,a]_\frkg)-\mu^{0}(\dM_\frkg[x,a]_\frkg)-\dM_{\frkh}\circ\mu_{1}(a)\circ\mu^{0}(x)+\dM_{\frkh}\circ\mu^{1}(x)\circ\mu_{1}(a)+\dM_{\frkh}\circ\mu_{1}([a,x]_\frkg)\\
&&-\mu^{0}(x)\circ\mu^{0}(\dM_\frkg a)+\mu^{0}(\dM_\frkg a)\circ\mu^{0}(x)+\mu^{0}([x,\dM_\frkg a]_\frkg)-\mu^{0}(x)\circ\dM_\frkh\circ\mu_{1}(a)+\mu^{0}(x)\circ\mu^{0}(\dM_\frkg a)\\
&&+\dM_\frkh\circ\mu_{1}(a)\circ\mu^{0}(x)-\mu^{0}(\dM_\frkg a)\circ\mu^{0}(x)\big)u\\
&=&0.
\end{eqnarray*}
Similarly, for all  $x\in\g_0,a\in\g_1,m\in\frkh_1$,   we have
\begin{eqnarray*}
 [\varphi([x,a]_\frkg)+\dM_\frkh(\nu(x,a))-\omega(x,\dM_\frkg a)-\mu^{0}(x)\varphi(a),m]_{\frkh}=0.
\end{eqnarray*}
 Therefore, we have
\begin{eqnarray}
\label{c1}\varphi([x,a]_\frkg)+\dM_\frkh(\nu(x,a))-\omega(x,\dM_\frkg a)-\mu^{0}(x)\varphi(a)\in\cen^0(\frkh).
\end{eqnarray}
Since $\cen(\frkh)=0$, we obtain \eqref{p10}. For all $a,b\in\g_1,u\in\frkh_0$, by\eqref{lift3} and \eqref{lift5}, we have
\begin{eqnarray*}
&&[\nu(\dM_\frkg a,b)-\mu_{1}(b)\varphi(a)-\nu(a,\dM_\frkg b)-\mu_{1}(a)\varphi(b),u]_\frkh\\
&=&[\nu(\dM_\frkg a,b),u]_{\frkh}-[\mu_{1}(b)\varphi(a),u]_{\frkh}-[\nu(a,\dM_\frkg b),u]_{\frkh}-[\mu_{1}(a)\varphi(b),u]_{\frkh}\\
&=&-\ad_1({\nu(b,\dM_\frkg a)})(u)-(\mu_{1}(b)[\varphi(a),u]_{\frkh}-[\varphi(a),\mu_{1}(b)u]_{\frkh})\\
&&-\ad_1({\nu(a,\dM_\frkg b)})(u)-
(\mu_{1}(a)[\varphi(b),u]_{\frkh}-[\varphi(b),\mu_{1}(a)u]_{\frkh})\\
&=&\big(-\mu_{1}(b)\circ\mu^{0}(\dM_\frkg a)+\mu^{1}(\dM_\frkg a)\circ\mu_{1}(b)+\mu_{1}([b,\dM_\frkg a]_\frkg)\\
&&-\mu_{1}(b)\circ\dM_\frkh\circ\mu_{1}(a)+\mu_{1}(b)\circ\mu^{0}(\dM_\frkg a)+\mu_{1}(a)\circ \dM_\frkh\circ\mu_{1}(b)-\mu^{1}(\dM_\frkg a)\circ\mu_{1}(b)\\
&&-\mu_{1}(a)\circ\mu^{0}(\dM_\frkg b)+\mu^{1}(\dM_\frkg b)\circ\mu_{1}(a)+\mu_{1}([a,\dM_\frkg b]_\frkg)\\
&&-\mu_{1}(a)\circ\dM_\frkh\circ\mu_{1}(b)+\mu_{1}(a)\circ\mu^{0}(\dM_\frkg b)+\mu_{1}(b)\circ \dM_\frkh\circ\mu_{1}(a)-\mu^{1}(\dM_\frkg b)\circ\mu_{1}(a)\big)u\\
&=&0.
\end{eqnarray*}
Thus, we get
\begin{eqnarray}
\label{c2}\nu(\dM_\frkg a,b)-\mu_{1}(b)\varphi(a)-\nu(a,\dM_\frkg b)-\mu_{1}(a)\varphi(b)\in\cen^1(\frkh),
\end{eqnarray}
which implies that \eqref{p11} holds.
For all $x,y,z\in\g_0,m\in\frkh_1$, by \eqref{lift4}, we have
\begin{eqnarray*}
&&[\mu^{0}(x)\omega(y,z)-\omega([x,y]_\frkg,z)+c.p.(x,y,z),m]_{\frkh}\\
&=&(\mu^{1}(x)[\omega(y,z),m]_{\frkh}-[\omega(y,z),\mu^{1}(x)(m)]_{\frkh})-\ad^1_{\omega([x,y]_\frkg,z)}(m)+c.p.(x,y,z)\\
&=&\big(\mu^{1}(x)\circ[\mu^{1}(y),\mu^{1}(z)]-\mu^{1}(x)\circ\mu^{1}([y,z]_\frkg)-[\mu^{1}(y),\mu^{1}(z)]\circ\mu^{1}(x)+\mu^{1}([y,z]_\frkg)\circ\mu^{1}(x)\\
&&-[\mu^{1}([x,y]_\g),\mu^{1}(z)]+\mu^{1}([[x,y]_\g,z]_\frkg)+c.p.(x,y,z)\big)m\\
&=&0.
\end{eqnarray*}
Similarly, for all $x,y,z\in\g_0,u\in\frkh_0$,  we have
$$[\mu^{0}(x)\omega(y,z)-\omega([x,y]_\frkg,z)+c.p.(x,y,z),u]_{\frkh}=0.$$
Therefore, we have
\begin{eqnarray}
\label{c3}\mu^{0}(x)\omega(y,z)-\omega([x,y]_\frkg,z)+c.p.(x,y,z)\in\cen^0(\frkh).
\end{eqnarray}
Since $\cen(\frkh)=0$, we obtain \eqref{p12}.
Similarly, for all $x,y\in\g_0,a\in\g_1,u\in\frkh_0$, by \eqref{lift4}-\eqref{lift5}, we have
\begin{eqnarray*}
~[\nu([x,y]_\frkg,a)+\nu([y,a]_\frkg,x)+\nu([a,x]_\frkg,y)-\mu^{1}(x)\nu(y,a)-\mu^{1}(y)\nu(a,x)-\mu_{1}(a)\omega(x,y),u]_{\frkh} =0.
\end{eqnarray*}
Therefore, we have
\begin{eqnarray}
\label{c4}\nu([x,y]_\frkg,a)+c.p.(x,y,a)-\mu_{1}(x)\nu(y,a)-\mu_{1}(y)\nu(a,x)-\mu_{2}(a)\omega(x,y)\in\cen^1(\frkh).
\end{eqnarray}
Since $\cen(\frkh)=0$, we obtain \eqref{p13}. Thus, we deduce that \eqref{p1}-\eqref{p13} hold. By Proposition \ref{extension},  $(\frkg\oplus\frkh, \dM_{\frkg\oplus\frkh}, [\cdot, \cdot]_{\frkg\oplus\frkh})$ is an extension of $\g$ by $\frkh$.

If we choose another section $s'$ of $\pi$, we obtain another extension $(\frkg\oplus\frkh, \dM_{\frkg\oplus\frkh}', [\cdot, \cdot]_{\frkg\oplus\frkh}')$  of $\g$ by $\frkh$. Obviously, we have
\begin{eqnarray*}
\pi_0(\mu_0'(x)-\mu_0(x))=\pi_0(s_0'(\bar{\mu}_0(x))-s_0(\bar{\mu}_0(x)))=\bar{\mu}_0(x)-\bar{\mu}_0(x)=0,
\end{eqnarray*}
which implies that $\mu_0'(x)-\mu_0(x)\in\SInn^0(\frkh)$. Since $\cen(\frkh)=0$, there is a unique linear map $\xi:\g_0\lon\frkh_0$ such that
\begin{eqnarray}\label{rep1}
\mu_0'(x)-\mu_0(x)=\ad_0({\xi(x)}).
\end{eqnarray}
Therefore, we obtain \eqref{iso1}. Similarly, there is a unique linear map $\eta:\g_1\lon\frkh_1$ such that
\begin{eqnarray}\label{rep2}
\mu_1'(a)-\mu_1(a)=\ad_1(\eta(a)).
\end{eqnarray}
Therefore, we obtain \eqref{iso2}. By Lemma \ref{adjoint} and \eqref{lift3}, for all $a\in\g_1$, we have
\begin{eqnarray*}
\ad_0(\varphi'(a)-\varphi(a))&=&\big(\delta(\mu_1'(a))-\mu_0'(\dM_\g a)\big)-\big(\delta(\mu_1(a))-\mu_0(\dM_\g a)\big)\\
                                    &=&\delta(\ad_1(\eta(a)))-\ad_0(\xi(\dM_\g a))\\
                                    &=&\ad_0(\dM_{\frkh}(\eta(a)))-\ad_0(\xi(\dM_\g a))\\
                                    &=&\ad_0(\dM_{\frkh}(\eta(a))-\xi(\dM_\g a)).
\end{eqnarray*}
Since $\cen(\frkh)=0$,  we obtain \eqref{iso3}. By Lemma \ref{ideal} and \eqref{lift4}, for all $x,y\in\g_0$ we have
\begin{eqnarray*}
&&\ad_0(\omega'(x,y)-\omega(x,y))\\&=&\big([\mu_0'(x),\mu_0'(y)]_C-\mu_0'([x,y]_\g)\big)-\big([\mu_0(x),\mu_0(y)]_C-\mu_0([x,y]_\g)\big)\\
                               &=&[\mu_0(x)+\ad_0(\xi(x)),\mu_0(y)+\ad_0(\xi(y))]_C-[\mu_0(x),\mu_0(y)]_C-\ad_0(\xi([x,y]_\g))\\
                               &=&[\mu_0(x),\ad_0(\xi(y))]_C+[\ad_0(\xi(x)),\mu_0(y)]_C+[\ad_0(\xi(x)),\ad_0(\xi(y))]_C-\ad_0(\xi([x,y]_\g))\\
                               &=&\ad_0(\mu^0(x)(\xi(y)))-\ad_0(\mu^0(y)(\xi(x)))+\ad_0([\xi(x),\xi(y)]_{\frkh})-\ad_0(\xi([x,y]_\g)),\\
                               &=&\ad_0(\mu^0(x)(\xi(y))-\mu^0(y)(\xi(x))+[\xi(x),\xi(y)]_{\frkh}-\xi([x,y]_\g)).
\end{eqnarray*}
By $\cen(\frkh)=0$, we obtain \eqref{iso4}. By Lemma \ref{ideal} and \eqref{lift5}, for all $x\in\g_0,a\in\g_1$, we have
\begin{eqnarray*}
&&\ad_1(\nu'(a,x)-\nu(a,x))\\&=&\big([\mu_1'(a),\mu_0'(x)]_C-\mu_1'([a,x]_\g)\big)-\big([\mu_1(a),\mu_0(x)]_C-\mu_1([a,x]_\g)\big)\\
                               &=&[\mu_1(a)+\ad_1(\eta(a)),\mu_0(x)+\ad_0(\xi(x))]_C-[\mu_1(a),\mu_0(x)]_C-\ad_1(\eta([a,x]_\g))\\
                               &=&[\mu_1(a),\ad_0(\xi(x))]_C+[\ad_1(\eta(a)),\mu_0(x)]_C+[\ad_1(\eta(a)),\ad_0(\xi(x))]_C-\ad_1(\eta([a,x]_\g))\\
                               &=&\ad_1(\mu_1(a)(\xi(x)))-\ad_1(\mu^1(x)(\eta(a)))+\ad_1([\eta(a),\xi(x)]_{\frkh})-\ad_1(\eta([a,x]_\g))\\
                               &=&\ad_1(\mu_1(a)(\xi(x))-\mu^1(x)(\eta(a))+[\eta(a),\xi(x)]_{\frkh}-\eta([a,x]_\g)).
\end{eqnarray*}
By $\cen(\frkh)=0$, we obtain \eqref{iso5}. Thus, we have \eqref{iso1}-\eqref{iso5}. Therefore, we deduce that
$(\frkg\oplus\frkh, \dM_{\frkg\oplus\frkh}, [\cdot, \cdot]_{\frkg\oplus\frkh})$ and $(\frkg\oplus\frkh, \dM_{\frkg\oplus\frkh}', [\cdot, \cdot]_{\frkg\oplus\frkh}')$ are isomorphic non-abelian extensions of $\g$ by $\frkh$. The proof is finished. \qed

\section{Obstruction of existence of non-abelian extensions of strict Lie 2-algebras}

In this section, given a homomorphism $\bar{\mu}:\g\lon\SOut(\h)$, where $\cen(\frkh)\not=0$, we consider the obstruction of existence of non-abelian extensions.
By choosing a section $s$ of $\pi:\SDer(\h)\lon\SOut(\h)$, we can still define $\mu_0(x)$ by \eqref{lift1}, and
$\mu_1(a)$ by \eqref{lift2}. Moreover, we can also choose linear maps $\varphi:\g_1\lon\frkh_0,~\omega:\wedge^2\mathfrak{g}_{0}\lon\frkh_0$
and $\nu:\mathfrak{g}_{0}\wedge\mathfrak{g}_{1}\lon\frkh_1$ such that \eqref{lift3}-\eqref{lift5} hold\footnote{Now $\varphi,~\omega,~\nu$ are not unique.}. Thus, $(\frkg\oplus\frkh, \dM_{\frkg\oplus\frkh}, [\cdot, \cdot]_{\frkg\oplus\frkh})$ is a non-abelian extension of $\g$ by $\frkh$ if and only if
\begin{eqnarray}\label{eq:extensiblecon}
\left\{\begin{array}{rcl}
\nu(\dM_\frkg a,b)-\mu_{1}(b)\varphi(a)-\nu(a,\dM_\frkg b)-\mu_{1}(a)\varphi(b)&=&0,\\
\varphi([x,a]_\frkg)+\dM_\frkh(\nu(x,a))-\omega(x,\dM_\frkg a)-\mu^{0}(x)\varphi(a)&=&0,\\
\mu^{0}(x)\omega(y,z)-\omega([x,y]_\frkg,z)+c.p.(x,y,z)&=&0,\\
\nu([x,y]_\frkg,a)+c.p.(x,y,a)-\mu^{1}(x)\nu(y,a)-\mu^{1}(y)\nu(a,x)-\mu_{1}(a)\omega(x,y)&=&0.
\end{array}\right.
\end{eqnarray}
 We denote by $\lambda=\varphi+\omega+\nu$. Let $D_{\mu}$ be the formal coboundary operator associated to  $\mu=(\mu_0,\mu_1)$. Then $\Omega=\D_{\mu}\lambda$ has four components as follows:
\begin{equation*}
\left\{
\begin{array}{lcll}
{\Omega_1} &=&\dM_{\rho}^{(1,0)}\varphi+\hat{\dM_\g}\omega+\hat{\partial}\nu   &\in\Hom(\mathfrak{g}_{0}\wedge\mathfrak{g}_{1},\frkh_0),\\
{\Omega_2} &=&\dM_{\rho}^{(0,1)}\varphi+\hat{\dM_\g}\nu                        &\in\Hom(\odot^2\g_1,\frkh_1),\\
{\Omega_3} &=&\dM_{\rho}^{(1,0)}\omega                                         &\in\Hom(\wedge^3\g_0,\frkh_0),\\
{\Omega_4} &=&\dM_{\rho}^{(0,1)}\omega+\dM_{\rho}^{(1,0)}\nu                   &\in\Hom(\wedge^2\g_0\wedge\g_1,\frkh_1).
\end{array}
\right.
\end{equation*}
More precisely, for all $x,y,z\in\g_0,a,b\in\g_1$, we have
\begin{equation*}
\left\{
\begin{array}{lcl}
{\Omega_1}(x,a) &=&\mu^0(x)\varphi(a)-\varphi([x,a]_\g)+\omega(x,\dM_\g a)-\dM_{\frkh}(\nu(x,a)),\\
{\Omega_2}(a,b) &=&\mu_1(a)\varphi(b)+\mu_1(b)\varphi(a)-\nu(\dM_\g a,b)-\nu(\dM_\g b,a),\\
{\Omega_3}(x,y,z) &=&\mu^0(x)\omega(y,z)-\omega([x,y]_\g,z)+c.p.(x,y,z),\\
{\Omega_4}(x,y,a) &=&\mu_1(a)\omega(x,y)+\mu^1(x)\nu(y,a)-\mu^1(y)\nu(x,a)-\big(\nu([x,y]_\g,a)+c.p.(x,y,a)\big).\\
\end{array}
\right.
\end{equation*}
Therefore,   $(\frkg\oplus\frkh, \dM_{\frkg\oplus\frkh}, [\cdot, \cdot]_{\frkg\oplus\frkh})$ is an extension of $\g$ by $\frkh$ if and only if $\Omega=\D_{\mu}\lambda=0.$

\begin{defi}
 Let $\bar{\mu}:\g\lon\SOut(\frkh)$ be a strict Lie $2$-algebra homomorphism. We call $\bar{\mu}$ an {\bf extensible homomorphism} if there exists a section $s$ of $\pi:\SDer(\h)\lon\SOut(\h)$ and linear maps $\varphi:\g_1\lon\frkh_0,~\omega:\wedge^2\mathfrak{g}_{0}\lon\frkh_0$
and $\nu:\mathfrak{g}_{0}\wedge\mathfrak{g}_{1}\lon\frkh_1$ such that \eqref{lift3}-\eqref{lift5} and \eqref{eq:extensiblecon} hold.
\end{defi}

For all $u\in\cen^0(\frkh),v\in\frkh_0,m\in\frkh_1$, we have
\begin{equation*}
[\mu_1(a)u,v]_{\frkh}=\mu_1(a)[u,v]_{\frkh}-[u,\mu_1(a)v]_{\frkh}=0.
\end{equation*}
Thus, we have $\mu_1(a)u\in\cen^1(\frkh)$. Moreover, we have
\begin{eqnarray*}
[\mu^0(x)u,v]_{\frkh}           &=&\mu^0(x)[u,v]_{\frkh}-[u,\mu^0(x)v]_{\frkh}=0,\\
\label{z}[\mu^0(x)u,m]_{\frkh}&=&\mu^1(x)[u,m]_{\frkh}-[u,\mu^1(x)m]_{\frkh}=0.
\end{eqnarray*}
Thus, we have $\mu^0(x)u\in\cen^0(\frkh)$. For $n\in\cen^1(\frkh)$, we have
\begin{eqnarray*}
[\mu^1(x)n,v]_{\frkh}=\mu^1(x)[n,v]_{\frkh}-[n,\mu^0(x)v]_{\frkh}=0.
\end{eqnarray*}
Thus, we have $\mu^1(x)n\in\cen^1(\frkh)$. Therefore, we can define $\hat{\mu}=(\hat{\mu}_0,\hat{\mu}_1):\g\lon\End(\cen(\frkh))$ by
\begin{eqnarray*}
&&\hat{\mu}_0(x)\triangleq\mu_0(x)|_{\cen(\frkh)}=\big(\mu^0(x)|_{\cen^0(\frkh)},\mu^1(x)|_{\cen^1(\frkh)}\big),\\
&&\hat{\mu}_1(a)\triangleq\mu_1(a)|_{\cen^0(\frkh)}.
\end{eqnarray*}
By \eqref{lift3}-\eqref{lift5},   $\hat{\mu}$ is a strict representation of $\g$ on the 2-term complex $\cen(\frkh)$. Moreover, by \eqref{rep1}-\eqref{rep2}, we deduce that different sections of $\pi$ give the same representation of $\g$ on $\cen(\frkh)$. In the sequel, we always assume that $\hat{\mu}$ is a representation of $\g$  on $\cen(\frkh)$,   which is induced by $\bar{\mu}$. By \eqref{c1}-\eqref{c4}, we have ${\Omega_1}(x,a),{\Omega_3}(x,y,z)\in\cen^0(\frkh)$ and ${\Omega_2}(a,b),{\Omega_4}(x,y,a)\in\cen^1(\frkh)$. Thus, we have $\Omega=\D_{\mu}\lambda\in C^3(\g;\cen(\frkh))$. Moreover, we have the following lemma.

\begin{lem}\label{obstruction}
$\Omega=\D_{\mu}\lambda$ is a $3$-cocycle on $\g$ with the coefficient in $\cen(\h)$ and the cohomology class $[\Omega]$ does not depend on the choices of the section $s$ of $\pi:\SDer(\h)\lon\SOut(\h)$ and  $\varphi,\omega,\nu$ that we made.
\end{lem}

\pf Denote by $\Upsilon=\D_{\hat{\mu}}\Omega$. It has five components as follows:
\begin{equation*}
\left\{
\begin{array}{lcll}
{\Upsilon_1} &=&\hat{\dM_\g}\Omega_1+\hat{\partial}\Omega_2                                  &\in\Hom(\odot^2\g_1,\cen^0(\frkh)),\\
{\Upsilon_2} &=&\dM_{\rho}^{(0,1)}\Omega_1+\dM_{\rho}^{(1,0)}\Omega_2+\hat{\dM_\g}\Omega_4 &\in\Hom(\g_0\wedge\odot^2\g_1,\cen^1(\frkh)),\\
{\Upsilon_3} &=&\dM_{\rho}^{(1,0)}\Omega_1+\hat{\dM_\g}\Omega_3+\hat{\partial}\Omega_4     &\in\Hom(\wedge^2\g_0\wedge\g_1,\cen^0(\frkh)),\\
{\Upsilon_4} &=&\dM_{\rho}^{(0,1)}\Omega_3+\dM_{\rho}^{(1,0)}\Omega_4                      &\in\Hom(\wedge^3\g_0\wedge\g_1,\cen^1(\frkh)),\\
{\Upsilon_5} &=&\dM_{\rho}^{(1,0)}\Omega_3                                                 &\in\Hom(\wedge^4\g_0,\cen^0(\frkh)).
\end{array}
\right.
\end{equation*}
For all $a,b\in\g_1$, by \eqref{lift3}, we have
\begin{eqnarray*}
\Upsilon_1(a,b)&=&(-1)^1\big(\Omega_1(\dM_\g a,b)+\Omega_1(\dM_\g b,a)\big)+\dM_{\frkh}\Omega_2(a,b)\\
               &=&-\mu^0(\dM_\g a)\varphi(b)+\varphi([\dM_\g a,b]_\g)-\omega(\dM_\g a,\dM_\g b)+\dM_{\frkh}(\nu(\dM_\g a,b))\\
               &&-\mu^0(\dM_\g b)\varphi(a)+\varphi([\dM_\g b,a]_\g)-\omega(\dM_\g b,\dM_\g a)+\dM_{\frkh}(\nu(\dM_\g b,a))\\
               &&+\dM_{\frkh}(\mu_1(a)\varphi(b))+\dM_{\frkh}(\mu_1(b)\varphi(a))-\dM_{\frkh}(\nu(\dM_\g a,b))-\dM_{\frkh}(\nu(\dM_\g b,a))\\
               &=&0.
\end{eqnarray*}
For all $x\in\g_0,a,b\in\g_1$, by \eqref{lift3} and \eqref{lift5}, we have
\begin{eqnarray*}
\Upsilon_2(x,a,b)&=&(-1)^1\big(\mu_1(a)\Omega_1(x,b)+\mu_1(b)\Omega_1(x,a)\big)+\big(\mu^1(x)\Omega_2(a,b)+(-1)^1\Omega_2([x,a]_\g,b)\\
                 &&+(-1)^1\Omega_2(a,[x,b]_\g)\big)+(-1)^2\big(\Omega_4(x,\dM_\g a,b)+\Omega_4(x,\dM_\g b,a)\big)\\
                 &=&-(\mu_1(a)\circ\mu^0(x))\varphi(b)+\mu_1(a)\varphi([x,b]_\g)-\mu_1(a)\omega(x,\dM_\g b)+(\mu_1(a)\circ\dM_{\frkh})\nu(x,b)\\
                 &&-(\mu_1(b)\circ\mu^0(x))\varphi(a)+\mu_1(b)\varphi([x,a]_\g)-\mu_1(b)\omega(x,\dM_\g a)+(\mu_1(b)\circ\dM_{\frkh})\nu(x,a)\\
                 &&+(\mu^1(x)\circ\mu_1(a))\varphi(b)+(\mu^1(x)\circ\mu_1(b))\varphi(a)-\mu^1(x)\nu(\dM_\g a,b)-\mu^1(x)\nu(\dM_\g b,a)\\
                 &&-\mu_1([x,a]_\g)\varphi(b)-\mu_1(b)\varphi([x,a]_\g)+\nu(\dM_\g[x,a]_\g,b)+\nu(\dM_\g b,[x,a]_\g)\\
                 &&-\mu_1(a)\varphi([x,b]_\g)-\mu_1([x,b]_\g)\varphi(a)+\nu(\dM_\g a,[x,b]_\g)+\nu(\dM_\g[x,b]_\g,a)\\
                 &&+\mu_1(b)\omega(x,\dM_\g a)+\mu^1(x)\nu(\dM_\g a,b)-\mu^1(\dM_\g a)\nu(x,b)-\nu([x,\dM_\g a]_\g,b)\\
                 &&-\nu(\dM_\g a,[x,b]_\g)+\nu(x,[\dM_\g a,b]_\g)+\mu_1(a)\omega(x,\dM_\g b)+\mu^1(x)\nu(\dM_\g b,a)\\
                 &&-\mu^1(\dM_\g b)\nu(x,a)-\nu([x,\dM_\g b]_\g,a)-\nu(\dM_\g b,[x,a]_\g)+\nu(x,[\dM_\g b,a]_\g)\\
                 &=&0.
\end{eqnarray*}
For all $x,y\in\g_0,a\in\g_1$, by \eqref{lift3} and \eqref{lift4}, we have
\begin{eqnarray*}
\Upsilon_3(x,y,a)&=&\mu^0(x)\Omega_1(y,a)-\mu^0(y)\Omega_1(x,a)-\Omega_1([x,y]_\g,a)-\Omega_1(y,[x,a]_\g)+\Omega_1(x,[y,a]_\g)\\
                 &&(-1)^3\Omega_3(x,y,\dM_\g a)+\dM_{\frkh}\Omega_4(x,y,a)\\
                 &=&(\mu^0(x)\circ\mu^0(y))\varphi(a)-\mu^0(x)\varphi([y,a]_\g)+\mu^0(x)\omega(y,\dM_\g a)-(\mu^0(x)\circ\dM_{\frkh})\nu(y,a)\\
                 &&-(\mu^0(y)\circ\mu^0(x))\varphi(a)+\mu^0(y)\varphi([x,a]_\g)-\mu^0(y)\omega(x,\dM_\g a)+(\mu^0(y)\circ\dM_{\frkh})\nu(x,a)\\
                 &&-\mu^0([x,y]_\g)\varphi(a)+\varphi([[x,y]_\g,a]_\g)-\omega([x,y]_\g,\dM_\g a)+\dM_{\frkh}\nu([x,y]_\g,a)\\
                 &&-\mu^0(y)\varphi([x,a]_\g)+\varphi([y,[x,a]_\g]_\g)-\omega(y,\dM_\g[x,a]_\g)+\dM_{\frkh}\nu(y,[x,a]_\g)\\
                 &&+\mu^0(x)\varphi([y,a]_\g)-\varphi([x,[y,a]_\g]_\g)+\omega(x,\dM_\g[y,a]_\g)-\dM_{\frkh}\nu(x,[y,a]_\g)\\
                 &&-\mu^0(x)\omega(y,\dM_\g a)-\mu^0(y)\omega(\dM_\g a,x)-\mu^0(\dM_\g a)\omega(x,y)+\omega([x,y]_\g,\dM_\g a)\\
                 &&+\omega([y,\dM_\g a]_\g,x)+\omega([\dM_\g a,x]_\g,y)+(\dM_{\frkh}\circ\mu_1(a))\omega(x,y)+(\dM_{\frkh}\circ\mu^1(x))\nu(y,a)\\
                 &&-(\dM_{\frkh}\circ\mu^1(y))\nu(x,a)-\dM_{\frkh}\nu([x,y]_\g,a)-\dM_{\frkh}\nu(y,[x,a]_\g)+\dM_{\frkh}\nu(x,[y,a]_\g)\\                  &=&0.
\end{eqnarray*}
For all $x,y,z\in\g_0,a\in\g_1$, by \eqref{lift4} and \eqref{lift5}, we have
\begin{eqnarray*}
&&\Upsilon_4(x,y,z,a)\\&=&(-1)^3\mu_1(a)\Omega_3(x,y,z)+\mu^1(x)\Omega_4(y,z,a)-\mu^1(y)\Omega_4(x,z,a)+\mu^1(z)\Omega_4(x,y,a)-\Omega_4([x,y]_\g,z,a)\\
                    &&+\Omega_4([x,z]_\g,y,a)-\Omega_4([y,z]_\g,x,a)-\Omega_4(y,z,[x,a]_\g)+\Omega_4(x,z,[y,a]_\g)-\Omega_4(x,y,[z,a]_\g)\\
                   &=&-\mu_1(a)\big(\mu^0(x)\omega(y,z)+\mu^0(y)\omega(z,x)+\mu^0(z)\omega(x,y)-\omega([x,y]_\g,z)-\omega([y,z]_\g,x)-\omega([z,x]_\g,y)\big)\\
                   &&+\mu^1(x)\big(\mu_1(a)\omega(y,z)+\mu^1(y)\nu(z,a)-\mu^1(z)\nu(y,a)-\nu([y,z]_\g,a)-\nu(z,[y,a]_\g)+\nu(y,[z,a]_\g)\big)\\
                   &&-\mu^1(y)\big(\mu_1(a)\omega(x,z)+\mu^1(x)\nu(z,a)-\mu^1(z)\nu(x,a)-\nu([x,z]_\g,a)-\nu(z,[x,a]_\g)+\nu(x,[z,a]_\g)\big)\\
                   &&+\mu^1(z)\big(\mu_1(a)\omega(x,y)+\mu^1(x)\nu(y,a)-\mu^1(y)\nu(x,a)-\nu([x,y]_\g,a)-\nu(y,[x,a]_\g)+\nu(x,[y,a]_\g)\big)\\
                   &&-\big(\mu_1(a)\omega([x,y]_\g,z)+\mu^1([x,y]_\g)\nu(z,a)-\mu^1(z)\nu([x,y]_\g,a)-\nu([[x,y]_\g,z]_\g,a)\\
                   &&-\nu(z,[[x,y]_\g,a]_\g)+\nu([x,y]_\g,[z,a]_\g)\big)+\big(\mu_1(a)\omega([x,z]_\g,y)+\mu^1([x,z]_\g)\nu(y,a)\\
                   &&-\mu^1(y)\nu([x,z]_\g,a)-\nu([[x,z]_\g,y]_\g,a)-\nu(y,[[x,z]_\g,a]_\g)+\nu([x,z]_\g,[y,a]_\g)\big)\\
                   &&-\big(\mu_1(a)\omega([y,z]_\g,x)+\mu^1([y,z]_\g)\nu(x,a)-\mu^1(x)\nu([y,z]_\g,a)-\nu([[y,z]_\g,x]_\g,a)\\
                   &&-\nu(x,[[y,z]_\g,a]_\g)+\nu([y,z]_\g,[x,a]_\g)\big)-\big(\mu_1([x,a]_\g)\omega(y,z)+\mu^1(y)\nu(z,[x,a]_\g)\\
                   &&-\mu^1(z)\nu(y,[x,a]_\g)-\nu([y,z]_\g,[x,a]_\g)-\nu(z,[y,[x,a]_\g]_\g)+\nu(y,[z,[x,a]_\g]_\g)\big)\\
                   &&+\big(\mu_1([y,a]_\g)\omega(x,z)+\mu^1(x)\nu(z,[y,a]_\g)-\mu^1(z)\nu(x,[y,a]_\g)-\nu([x,z]_\g,[y,a]_\g)\\
                   &&-\nu(z,[x,[y,a]_\g]_\g)+\nu(x,[z,[y,a]_\g]_\g)\big)-\big(\mu_1([z,a]_\g)\omega(x,y)+\mu^1(x)\nu(y,[z,a]_\g)\\
                   &&-\mu^1(y)\nu(x,[z,a]_\g)-\nu([x,y]_\g,[z,a]_\g)-\nu(y,[x,[z,a]_\g]_\g)+\nu(x,[y,[z,a]_\g]_\g)\big)\\
                   &=&0.
\end{eqnarray*}
 For all $x,y,z,w\in\g_0$, by the  similar computation, we have $\Upsilon_5(x,y,z,w)=0$. Thus, we obtain $\Upsilon=\D_{\hat{\mu}}\Omega=0$.

Now Let us check  that the cohomology class $[\Omega]=[\D_{\mu}(\varphi+\omega+\nu)]$ does not depend on the choices of the section $s$ of $\pi:\SDer(\h)\lon\SOut(\h)$ and  $\varphi,\omega,\nu$ that we made. Let $s'$ be another section of $\pi$,  we have $\mu'=(\mu'_0,\mu'_1)=(s_0'\circ\bar{\mu}_0,s_1'\circ\bar{\mu}_1)$, where $\mu'_0(x)=(\mu'^0(x),\mu'^1(x))\in\SDer^0(\frkg),\mu'_1(a)\in\SDer^1(\frkg)$ and choose $\varphi',\omega',\nu'$ such that
\eqref{lift3}-\eqref{lift5} hold. Let $\Omega'=\D_{\mu'}(\varphi'+\omega'+\nu')$. We are going to prove that $[\Omega]=[\Omega']$.
Since $s$ and $s'$ are sections of $\pi$, we have linear maps $\xi:\g_0\lon\frkh_0$ and $\eta:\g_1\lon\frkh_1$ such that
\begin{eqnarray*}
\mu'_0(x)-\mu_0(x)=\ad_0({\xi(x)}),\quad
\mu'_1(a)-\mu_1(a)=\ad_1({\eta(a)}).
\end{eqnarray*}
We define $\varphi^*,\omega^*,\nu^*$ by
\begin{eqnarray*}
\varphi^*(a)&=&\varphi(a)+\dM_{\frkh}\eta(a)-\xi(\dM_\g a),\\
\omega^*(x,y)&=&\omega(x,y)+\mu^0(x)\xi(y)-\mu^0(y)\xi(x)-\xi([x,y]_\g)+[\xi(x),\xi(y)]_{\frkh},\\
\nu^*(a,x)&=&\nu(a,x)+\mu_1(a)\xi(x)-\mu^1(x)\eta(a)-\eta([a,x]_\g)+[\eta(a),\xi(x)]_{\frkh}.
\end{eqnarray*}
By straightforward computations, we obtain that \eqref{lift3}-\eqref{lift5} hold for $\mu',\varphi^*,\omega^*,\nu^*$. Let
$$\Omega^*=\D_{\mu'}(\varphi^*+\omega^*+\nu^*)=\Omega^*_1+\Omega^*_2+\Omega^*_3+\Omega^*_4.$$
For all $x\in\g_0,a\in\g_1$, by \eqref{lift3}, we have
\begin{eqnarray*}
&&\Omega^*_1(x,a)\\
&=&\mu'^0(x)\varphi^*(a)-\varphi^*([x,a]_\g)+\omega^*(x,\dM_\g a)-\dM_{\frkh}\nu^*(x,a)\\
\emptycomment{
              &=&\big(\mu^0(x)+\ad^0({\xi(x)})\big)\big(\varphi(a)+\dM_{\frkh}\eta(a)-\xi(\dM_\g a)\big)-\big(\varphi([x,a]_\g)+\dM_{\frkh}\eta([x,a]_\g)-\xi(\dM_\g[x,a]_\g)\big)\\
              &&+\big(\omega(x,\dM_\g a)+\mu^0(x)\xi(\dM_\g a)-\mu^0(\dM_\g a)\xi(x)-\xi([x,\dM_\g a]_\g)+[\xi(x),\xi(\dM_\g a)]_{\frkh}\big)\\
              &&+\dM_{\frkh}\big(\nu(a,x)+\mu_1(a)\xi(x)-\mu^1(x)\eta(a)-\eta([a,x]_\g)+[\eta(a),\xi(x)]_{\frkh}\big)\\
              }
              &=&\mu^0(x)\varphi(a)+(\mu^0(x)\circ\dM_{\frkh})\eta(a)-\mu^0(x)\xi(\dM_\g a)+[\xi(x),\varphi(a)]_{\frkh}+[\xi(x),\dM_{\frkh}\eta(a)]_{\frkh}-[\xi(x),\xi(\dM_\g a)]_{\frkh}\\
              &&-\varphi([x,a]_\g)-\dM_{\frkh}\eta([x,a]_\g)+\xi(\dM_\g[x,a]_\g)+\omega(x,\dM_\g a)+\mu^0(x)\xi(\dM_\g a)-\mu^0(\dM_\g a)\xi(x)\\
              &&-\xi([x,\dM_\g a]_\g)+[\xi(x),\xi(\dM_\g a)]_{\frkh}+\dM_{\frkh}\nu(a,x)+(\dM_{\frkh}\circ\mu_1(a))\xi(x)-(\dM_{\frkh}\circ\mu^1(x))\eta(a)\\
              &&-\dM_{\frkh}\eta([a,x]_\g)+\dM_{\frkh}[\eta(a),\xi(x)]_{\frkh}\\
               &=&\Omega_1(x,a).
\end{eqnarray*}
For all $a,b\in\g_1$, by \eqref{lift3}, we have
\begin{eqnarray*}
&&\Omega^*_2(a,b)\\&=&\mu_1'(a)\varphi^*(b)+\mu_1'(b)\varphi^*(a)-\nu^*(\dM_\g a,b)-\nu^*(\dM_\g b,a)\\
\emptycomment{
              &=&\big(\mu_1(a)+\ad_{\eta(a)}\big)\big(\varphi(b)+\dM_{\frkh}(\eta(b))-\xi(\dM_\g(b))\big)+\big(\mu_1(b)+\ad_{\eta(b)}\big)\big(\varphi(a)+\dM_{\frkh}(\eta(a))-\xi(\dM_\g(a))\big)\\
              &&+\big(\nu(b,\dM_\g(a))+\mu_1(b)(\xi(\dM_\g(a)))-\mu^1(\dM_\g(a))(\eta(b))-\eta([b,\dM_\g(a)]_\g)+[\eta(b),\xi(\dM_\g(a))]_{\frkh}\big)\\
              &&+\big(\nu(a,\dM_\g(b))+\mu_1(a)(\xi(\dM_\g(b)))-\mu^1(\dM_\g(b))(\eta(a))-\eta([a,\dM_\g(b)]_\g)+[\eta(a),\xi(\dM_\g(b))]_{\frkh}\big)\\
              }
              &=&\mu_1(a)\varphi(b)+(\mu_1(a)\circ\dM_{\frkh})\eta(b)-\mu_1(a)\xi(\dM_\g b)+[\eta(a),\varphi(b)]_{\frkh}+[\eta(a),\dM_{\frkh}\eta(b)]_{\frkh}-[\eta(a),\xi(\dM_\g b)]_{\frkh}\\
              &&+\mu_1(b)\varphi(a)+(\mu_1(b)\circ\dM_{\frkh})\eta(a)-\mu_1(b)\xi(\dM_\g a)+[\eta(b),\varphi(a)]_{\frkh}+[\eta(b),\dM_{\frkh}\eta(a)]_{\frkh}-[\eta(b),\xi(\dM_\g a)]_{\frkh}\\
              &&+\nu(b,\dM_\g a)+\mu_1(b)\xi(\dM_\g a)-\mu^1(\dM_\g a)\eta(b)-\eta([b,\dM_\g a]_\g)+[\eta(b),\xi(\dM_\g a)]_{\frkh}\\
              &&+\nu(a,\dM_\g b)+\mu_1(a)\xi(\dM_\g b)-\mu^1(\dM_\g b)\eta(a)-\eta([a,\dM_\g b]_\g)+[\eta(a),\xi(\dM_\g b)]_{\frkh}\\
              &=&\Omega_2(a,b).
\end{eqnarray*}
For all $x,y\in\g_0,a\in\g_1$, by \eqref{lift4} and \eqref{lift5}, we have
\begin{eqnarray*}
&&\Omega^*_4(x,y,a)\\&=&\mu_1'(a)\omega^*(x,y)+\mu'^1(x)\nu^*(y,a)-\mu'^1(y)\nu^*(x,a)-\nu^*([x,y]_\g,a)-\nu^*(y,[x,a]_\g)+\nu^*(x,[y,a]_\g)\\
                 &=&\big(\mu_1(a)+\ad_1({\eta(a)})\big)\big(\omega(x,y)+\mu^0(x)\xi(y)-\mu^0(y)\xi(x)-\xi([x,y]_\g)+[\xi(x),\xi(y)]_{\frkh}\big)\\
                 &&-\big(\mu^1(x)+\ad^1_{\xi(x)}\big)\big(\nu(a,y)+\mu_1(a)\xi(y)-\mu^1(y)\eta(a)-\eta([a,y]_\g)+[\eta(a),\xi(y)]_{\frkh}\big)\\
                 &&+\big(\mu^1(y)+\ad^1_{\xi(y)}\big)\big(\nu(a,x)+\mu_1(a)\xi(x)-\mu^1(x)\eta(a)-\eta([a,x]_\g)+[\eta(a),\xi(x)]_{\frkh}\big)\\
                 &&+\big(\nu(a,[x,y]_\g)+\mu_1(a)\xi([x,y]_\g)-\mu^1([x,y]_\g)\eta(a)-\eta([a,[x,y]_\g]_\g)+[\eta(a),\xi([x,y]_\g)]_{\frkh}\big)\\
                 &&+\big(\nu([x,a]_\g,y)+\mu_1([x,a]_\g)\xi(y)-\mu^1(y)\eta([x,a]_\g)-\eta([[x,a]_\g,y]_\g)+[\eta([x,a]_\g),\xi(y)]_{\frkh}\big)\\
                 &&-\big(\nu([y,a]_\g,x)+\mu_1([y,a]_\g)\xi(x)-\mu^1(x)\eta([y,a]_\g)-\eta([[y,a]_\g,x]_\g)+[\eta([y,a]_\g),\xi(x)]_{\frkh}\big)\\
                 &=&\Omega_4(x,y,a).
\end{eqnarray*}
Similarly, we have $\Omega^*_3(x,y,z)=\Omega_3(x,y,z).$ Therefore, we obtain $\Omega=\Omega^*.$ Since
the equations \eqref{lift3}-\eqref{lift5} hold for $\mu',\varphi^*,\omega^*,\nu^*$ and $\mu',\varphi',\omega',\nu'$ respectively, we have
\begin{eqnarray*}
\delta(\mu_1'(a))-\mu_0'(\dM_\g a)&=&\ad_0({\varphi'(a)})=\ad_0({\varphi^*(a)}),\\
\label{w}[\mu_0'(x),\mu_0'(y)]_C-\mu_0'([x,y]_\g)&=&\ad_0({\omega'(x,y)})=\ad_0({\omega^*(x,y)}),\\
\label{v}[\mu_1'(a),\mu_0'(x)]_C-\mu_1'([a,x]_\g)&=&\ad_1({\nu'(a,x)})=\ad_1({\nu^*(a,x)}).
\end{eqnarray*}
Thus, we can define
\begin{equation*}
\left\{
\begin{array}{lcll}
{t_1} &=&\varphi'-\varphi^* &\in\Hom(\g_1,\cen^0(\frkh)),\\
{t_2} &=&\omega'-\omega^*   &\in\Hom(\wedge^2\g_0,\cen^0(\frkh)),\\
{t_3} &=&\nu'-\nu^*         &\in\Hom(\g_0\wedge\g_1,\cen^1(\frkh)).\\
\end{array}
\right.
\end{equation*}
Denote by $T=t_1+t_2+t_3\in C^2(\g;\cen(\frkh))$, and we have
\begin{eqnarray*}
\Omega'=\D_{\mu'}(\varphi'+\omega'+\nu')&=&\D_{\mu'}(\varphi^*+\omega^*+\nu^*+T)\\
                                 &=&\D_{\mu'}(\varphi^*+\omega^*+\nu^*)+\D_{\mu'}(T)\\
                                 &=&\D_{\mu}(\varphi+\omega+\nu)+\D_{\hat{\mu}}(T)\\
                                 &=&\Omega+\D_{\hat{\mu}}(T).
\end{eqnarray*}
Therefore, we have $[\Omega]=[\Omega']$. The proof is finished. \qed\vspace{3mm}

Now we are ready to give the main result in this paper, namely the obstruction of a strict Lie $2$-algebra homomorphism $\bar{\mu}:\g\lon\SOut(\frkh)$ being extensible is given by the cohomology class $[\Omega]\in\huaH^3(\g;\cen(\frkh))_{\hat{\mu}}$.

\begin{thm}
Let $\bar{\mu}:\g\lon\SOut(\frkh)$ be a strict Lie $2$-algebra homomorphism. Then $\bar{\mu}$ is an extensible homomorphism if and only if $$[\Omega]=[\D_{\mu}(\varphi+\omega+\nu)]=[0].$$
\end{thm}

\pf Let $\bar{\mu}:\g\lon\SOut(\frkh)$ be an extensible strict Lie $2$-algebra homomorphism. Then we can choose a section $s$ of $\pi:\SDer(\h)\lon\SOut(\h)$ and define $\mu_0(x)=(\mu^{0}(x),\mu^{1}(x))$ by \eqref{lift1},
$\mu_1(a)$ by \eqref{lift2} respectively. Moreover, we can choose linear maps $\varphi:\g_1\lon\frkh_0,~\omega:\wedge^2\mathfrak{g}_{0}\lon\frkh_0$
and $\nu:\mathfrak{g}_{0}\wedge\mathfrak{g}_{1}\lon\frkh_1$ such that \eqref{lift3}-\eqref{lift5} hold. Since $\bar{\mu}$ is extensible, we have $\D_{\mu}(\varphi+\omega+\nu)=0$, which implies that $[\Omega]=[0].$

Conversely, if $[\D_{\mu}(\varphi+\omega+\nu)]=[0]$,
then there exists $\sigma=\sigma_1+\sigma_2+\sigma_3\in C^2(\g;\cen(\frkh))$, where $\sigma_1\in\Hom(\g_1,\cen^0(\frkh)),\sigma_2\in\Hom(\wedge^2\g_0,\cen^0(\frkh)),\sigma_3\in\Hom(\g_0\wedge\g_1,\cen^1(\frkh))$ such that $\D_{\hat{\mu}}(\sigma)=\D_{\mu}(\varphi+\omega+\nu)$. Thus, we have
\begin{eqnarray*}
\D_{\mu}(\varphi+\omega+\nu-\sigma)&=&\D_{\mu}(\varphi+\omega+\nu)-\D_{\mu}(\sigma)\\
                               &=&\D_{\mu}(\varphi+\omega+\nu)-\D_{\hat{\mu}}(\sigma)\\
                               &=&0.
\end{eqnarray*}
Since $\sigma=\sigma_1+\sigma_2+\sigma_3\in C^2(\g;\cen(\frkh))$, we also have
\begin{eqnarray*}
\delta(\mu_1(a))-\mu_0(\dM_\g a)&=&\ad_0({\varphi(a)-\sigma_1(a)}),\\
\label{u}[\mu_0(x),\mu_0(y)]_C-\mu_0([x,y]_\g)&=&\ad_0({\omega(x,y)-\sigma_2(x,y)}),\\
\label{t}[\mu_1(a),\mu_0(x)]_C-\mu_1([a,x]_\g)&=&\ad_1({\nu(a,x)-\sigma_3(a,x)}).
\end{eqnarray*}
By Proposition \ref{extension}, we can construct a strict Lie 2-algebra $(\frkg\oplus\frkh, \dM_{\frkg\oplus\frkh}, [\cdot, \cdot]_{\frkg\oplus\frkh})$ by $\mu,\varphi(a)-\sigma_1(a),\omega(x,y)-\sigma_2(x,y),\nu(a,x)-\sigma_3(a,x)$.
Therefore, $\bar{\mu}$ is an extensible homomorphism. The proof is finished. \qed\vspace{3mm}

The following theorem classifies non-abelian extensions of $\g$ by $\h$ once they exist.

\begin{thm}
Let $\bar{\mu}:\g\lon\SOut(\frkh)$ be an extensible homomorphism. Then isomorphism classes of non-abelian extensions of $\g$ by $\frkh$ induced by $\bar{\mu}$ are parameterized by $\huaH^2(\g;\cen(\frkh))_{\hat{\mu}}$.
\end{thm}

\pf Since $\bar{\mu}$ is an extensible homomorphism. We can choose a section $s$ of $\pi$ and define $\mu_0$, $\mu_1$ by \eqref{lift1}, \eqref{lift2} respectively. We choose linear maps $\varphi,\omega,\nu$
such that \eqref{lift3}-\eqref{lift5} hold and $\D_{\mu}(\varphi+\omega+\nu)=0$. Thus, the strict Lie 2-algebra $(\frkg\oplus\frkh, \dM_{\frkg\oplus\frkh}, [\cdot, \cdot]_{\frkg\oplus\frkh})$ define by \eqref{d}, \eqref{bo} and \eqref{b1} is a non-abelian extension of $\g$ by $\frkh$, which is induced by $\bar{\mu}$. Let $s'$ be another section of $\pi$ and define $\mu_0'$, $\mu_1'$ by \eqref{lift1}, \eqref{lift2}. We also choose linear maps $\varphi',\omega',\nu'$ such that \eqref{lift3}-\eqref{lift5} hold and $\D_{\mu'}(\varphi'+\omega'+\nu')=0$. Since $s$ and $s'$ are sections of $\pi$, we have linear maps $\xi:\g_0\lon\frkh_0$ and $\eta:\g_1\lon\frkh_1$ such that
\begin{eqnarray*}
\mu_0(x)=\mu_0'(x)+\ad_0({\xi(x)}),\quad
\mu_1(a)=\mu_1'(a)+\ad_1({\eta(a)}).
\end{eqnarray*}
We define $\varphi^*,\omega^*,\nu^*$ by
\begin{eqnarray*}
\varphi^*(a)&=&\varphi'(a)+\dM_{\frkh}\eta(a)-\xi(\dM_\g a),\\
\omega^*(x,y)&=&\omega'(x,y)+\mu'^0(x)\xi(y)-\mu'^0(y)\xi(x)-\xi([x,y]_\g)+[\xi(x),\xi(y)]_{\frkh},\\
\nu^*(a,x)&=&\nu'(a,x)+\mu_1'(a)(\xi(x))-\mu'^1(x)\eta(a)-\eta([a,x]_\g)+[\eta(a),\xi(x)]_{\frkh}.
\end{eqnarray*}
By the computation in Lemma \ref{obstruction}, we have $\D_{\mu}(\varphi^*+\omega^*+\nu^*)=\D_{\mu'}(\varphi'+\omega'+\nu')=0$. Thus, the strict Lie $2$-algebra $(\frkg\oplus\frkh, \dM_{\frkg\oplus\frkh}, [\cdot, \cdot]_{\frkg\oplus\frkh})$  constructed from $\mu,\varphi^*,\omega^*,\nu^*$ is isomorphic to the strict Lie $2$-algebra $(\frkg\oplus\frkh, \dM_{\frkg\oplus\frkh}, [\cdot, \cdot]_{\frkg\oplus\frkh})$  constructed from $\mu',\varphi',\omega',\nu'$. Thus, we only need to study the strict Lie $2$-algebras constructed from a fix section $s$. For all $\tilde{\varphi},\tilde{\omega},\tilde{\nu}$, which satisfy \eqref{lift3}-\eqref{lift5} and $\D_{\mu}(\tilde{\varphi}+\tilde{\omega}+\tilde{\nu})=0$. Thus, we can define
\begin{equation*}
\left\{
\begin{array}{lcll}
{\tau_1} &=&\varphi-\tilde{\varphi}       &\in\Hom(\g_1,\cen^0(\frkh)),\\
{\tau_2} &=&\omega-\tilde{\omega}         &\in\Hom(\wedge^2\g_0,\cen^0(\frkh)),\\
{\tau_3} &=&\nu-\tilde{\nu}               &\in\Hom(\g_0\wedge\g_1,\cen^1(\frkh)).
\end{array}
\right.
\end{equation*}
Denote by $\tau=\tau_1+\tau_2+\tau_3\in C^2(\g;\cen(\frkh))$, and we have
\begin{eqnarray*}
\D_{\hat{\mu}}(\tau)=\D_{\mu}\big((\varphi+\omega+\nu)-(\tilde{\varphi}+\tilde{\omega}+\tilde{\nu})\big)=0-0=0,
\end{eqnarray*}
which implies  that $\big((\varphi+\omega+\nu)-(\tilde{\varphi}+\tilde{\omega}+\tilde{\nu})\big)\in\huaZ^2(\g;\cen(\frkh))$.

Moreover, if the strict Lie $2$-algebra $(\frkg\oplus\frkh, \dM_{\frkg\oplus\frkh}, [\cdot, \cdot]_{\frkg\oplus\frkh})$  constructed from $\mu,\varphi,\omega,\nu$ is isomorphic to the strict Lie $2$-algebra $(\frkg\oplus\frkh, \dM_{\frkg\oplus\frkh}, [\cdot, \cdot]_{\frkg\oplus\frkh})$  constructed from $\mu,\tilde{\varphi},\tilde{\omega},\tilde{\nu}$. Then there exist linear maps $\xi:\g_0\lon\frkh_0$ and $\eta:\g_1\lon\frkh_1$ which does not change $\mu_0$ and $\mu_1$, i.e. $\xi:\g_0\lon\cen^0(\frkh)$ and $\eta:\g_1\lon\cen^1(\frkh)$, such that
\begin{eqnarray*}
\varphi(a)- \tilde{\varphi}(a)&=&\dM_{\frkh}\eta(a)-\xi(\dM_\g a),\\
{\omega}(x,y)-\tilde{\omega}(x,y)&=&\mu^0(x)\xi(y)-\mu^0(y)\xi(x)-\xi([x,y]_\g),\\
{\nu}(a,x)-\tilde{\nu}(a,x)&=&\mu_1(a)\xi(x)-\mu^1(x)\eta(a)-\eta([a,x]_\g),
\end{eqnarray*}
 which is equivalent to that
\begin{eqnarray*}
\big(({\varphi}+{\omega}+{\nu})-(\tilde{\varphi}+\tilde{\omega}+\tilde{\nu})\big)=\D_{\hat{\mu}}(\xi+\eta)\in\huaB^2(\g;\cen(\frkh)).
\end{eqnarray*}
Thus, isomorphism classes of non-abelian extensions of $\g$ by $\frkh$ induced by $\bar{\mu}$ are parameterized by $\huaH^2(\g;\cen(\frkh))_{\hat{\mu}}$.
 The proof is finished. \qed

\begin{cor}
The isomorphism classes of non-abelian extensions of a strict Lie $2$-algebra $\g$ by a strict Lie $2$-algebra $\frkh$ correspond bijectively to the set of pairs $(\bar{\mu},[\kappa])$, where $\bar{\mu}$ is an extensible homomorphism from $\g$ to $\SOut(\frkh)$ and $[\kappa]\in\huaH^2(\g;\cen(\frkh))_{\hat{\mu}}$.
\end{cor}

\end{document}